\documentclass[12pt,a4paper,reqno]{amsart}

\usepackage{amsmath,amsthm,amsfonts,amssymb,euscript,cite, mathrsfs,mathtools}
\usepackage{graphics}
\usepackage{xcolor} %
\usepackage{hyperref}
\hypersetup{colorlinks,linkcolor={blue},citecolor={blue},urlcolor={red}}  
\usepackage{enumitem} %
\usepackage{ulem}
\usepackage[noabbrev,capitalize,nameinlink]{cleveref}

\crefname{equation}{}{}
\crefname{algocf}{Algorithm}{Algorithms}
\crefname{equation}{}{} %
\crefname{algocf}{Algorithm}{Algorithms}

\usepackage{graphicx}
\usepackage{color}
\usepackage{transparent}

\usepackage{fullpage} %

\usepackage{seqsplit}

\definecolor{green}{rgb}{0,0.8,0} %

\definecolor{babypink}{rgb}{0.96,0.76,0.76}

\newcommand{\rt}{{\mathbb R^3}}

\newcommand{\pa}{\partial}

\newcommand{\grad}{\nabla_x}

\newcommand{\gab}{g^{\alpha\beta}}
\newcommand{\hab}{h^{\alpha\beta}}

\newcommand{\pab}{\partial_\beta}
\newcommand{\paa}{\partial_\alpha}

\newcommand{\pat}{\partial_t}

\newcommand{\pai}{\partial_i}

\newcommand{\pao}{\pa_\omega}

\newcommand{\ti}{\tilde}

\newcommand{\la}{\langle}
\newcommand{\ra}{\rangle}

\newcommand{\ls}{\lesssim}

\newcommand{\inv}{^{-1}}

\newcommand{\de}{\nabla_{t,x}} %

\newcommand{\f}{\frac}

\newcommand{\iy}{\infty}

\newcommand{\wm}{w_{(m)}}

\newcommand{\cinte}{\chi_\text{inte}}

\newcommand{\cone}{\chi^\text{cone}}
\newcommand{\cexte}{\chi_\text{exte}}

\newcommand{\crt}{{C^{ R}_{ T}}}
\newcommand{\cut}{{C^{ U}_{ T}}}
\newcommand{\cutt}{{\ti C^{ U}_{ T}}}
\newcommand{\crtt}{{\ti C^{ R}_{ T}}}
\newcommand{\ltcrt}{{L^2\lr{C^{ R}_{ T}}}}
\newcommand{\ltctr}{{L^2\lr{C^T_R}}}
\newcommand{\ltctrt}{{L^2\lr{\ti C^T_R}}}
\newcommand{\ltcrtt}{{L^2\lr{\crtt}}}
\newcommand{\ltcut}{{L^2\lr{C^{ U}_{ T}}}}

\newcommand{\ltcutt}{{L^2\lr{\cutt}}}
\newcommand{\licut}{{L^\iy\lr{C^{ U}_{ T}}}}

\newcommand{\licrt}{{L^\iy\lr{C^R_{ T}}}}

\newcommand{\lictr}{{L^\iy\lr{C^T_R}}}

\newcommand{\inte}{{C^{<3T/4}_{ T}}}
\newcommand{\tinte}{{\ti C^{<3T/4}_T}}

\newcommand{\rei}{{\calR_\text{1}}}
\newcommand{\reii}{{\calR_\text{2}}}
\newcommand{\co}{{D_{tr}^{R}}}
\newcommand{\dtr}{{D_{tr}}}

\newcommand{\lolt}{{L^1L^{2}}}
\newcommand{\ltlt}{{L^2L^2}}

\newcommand{\lt}{{L^2}}
\newcommand{\p}{\phi}

\newcommand{\supp}{\text{supp\,}}

\renewcommand{\pm}{\phi_{\le m}}
\newcommand{\pmo}{\phi_{\le m+1}}
\newcommand{\pmt}{\phi_{\le m+2}}
\newcommand{\pmn}{\p_{\le m+n}}
\newcommand{\nm}{\la t-r\ra} %
\newcommand{\mn}{\la r-t \ra} %
\newcommand{\on}{\la s-\rho\ra} %
\newcommand{\no}{\la \rho - s \ra} %
\newcommand{\jr}{\la r\ra} %
\newcommand{\ju}{\la u\ra}
\newcommand{\jxu}{\la\xu\ra}
\newcommand{\jt}{\la t\ra}
\newcommand{\jv}{\la v\ra}
\newcommand{\js}{\la s\ra}
\newcommand{\jrho}{\la\rho\ra}

\newcommand{\lr}[1]{\left( #1 \right)}

\let\arXiv\arxiv
\def\doi#1{ {\href{http://dx.doi.org/#1}
   {{\mdseries\ttfamily DOI}}}}

\newtheorem{theorem}{Theorem}[section]
\newtheorem{corollary}[theorem]{Corollary}
\newtheorem{lemma}[theorem]{Lemma}
\newtheorem{proposition}[theorem]{Proposition}
\theoremstyle{definition}
\crefname{claim}{Claim}{Claims}

\newtheorem{definition}[theorem]{Definition}

\theoremstyle{remark}
\newtheorem{remark}[theorem]{Remark}
\theoremstyle{conjecture}

\numberwithin{equation}{section}

\newcommand{\x}{\alpha}
\newcommand{\xb}{\beta}
\newcommand{\xd}{\delta}
	
\newcommand{\xg}{\gamma}
	
\newcommand{\eps}{\epsilon}

\newcommand{\xk}{\kappa}
\newcommand{\xl}{\lambda}
	
\newcommand{\xo}{\omega}
	\newcommand{\xO}{\Omega}

\newcommand{\xs}{\sigma}
	\newcommand{\xS}{\Sigma}
\newcommand{\xt}{\theta}
	
\newcommand{\xu}{\upsilon}

\newcommand{\N}{{\mathbb N}}

\newcommand{\R}{\mathbb R}

\newcommand{\bbS}{\mathbb S}

\newcommand{\calF}{\mathcal F}

\newcommand{\calR}{\mathcal R}

\newcommand{\bo}{\Box}

\usepackage{microtype} %

\begin{document}

\title{Pointwise decay for the wave equation on nonstationary spacetimes}

\author{Shi-Zhuo Looi}
\address{Department of Mathematics, University of Kentucky, Lexington, 
  KY  40506}
\email{Shizhuo.Looi@uky.edu}

\begin{abstract}
This first article in a two-part series (the second article being [arXiv:2205.13197]) assumes a weak local energy decay estimate holds and proves that solutions to the linear wave equation with variable coefficients in $\R^{1+3}$, first-order terms, and a potential decay at a rate depending on how rapidly the vector fields of the metric, first-order terms, and potential decay at spatial infinity. We prove results for both stationary and nonstationary metrics. The proof uses local energy decay to prove an initial decay rate, and then uses the one-dimensional reduction repeatedly to achieve the full decay rate.
\end{abstract}

\maketitle

\section{Introduction}\label{sec:intro}

In this paper, we examine pointwise decay for linear wave equations on asymptotically flat, nonstationary and stationary backgrounds in $1+3$ dimensions and show how, given certain \textit{weak local energy decay} estimates, the decay rate of the solution depends on the relative rates of the radial decay of the potential, the first-order coefficients and the background geometry. More on this assumption can be found in \cref{definition:weak}. 

Let \begin{equation}\label{P def}
P := \paa\gab(t,x)\pab + g^\xo(t,x) \Delta_\xo + \paa A^\x(t,x) + B^\x(t,x)\paa + V(t,x)
\end{equation} where the conditions on the potential $V$, the coefficients $A,B,g^\xo$ and the Lorentzian metric $g$ are given in the main result, \cref{thm:main}. $\Delta_\xo$ denotes the Laplace operator on the unit sphere.%
We let $\x,\xb$ range across $0, \dots, 3$. We consider the linear Cauchy problem
\begin{equation} \label{eq:problem}
P\p= f, \ \ (\phi(0), \bar N\phi(0)) =(\p_0,\p_1) 
\end{equation}where $\bar N$ denotes the unit normal derivative to the hypersurface $\{t=0\}$.

In the next three definitions, we consider a few types of bounds for \eqref{eq:problem} as a precursor to the main result, \cref{thm:main}. Given $x\in \R^3$, let $r :=|x|$ and $\jr := (1 + r^2)^{1/2}$. 

We will use the following norms throughout the paper. In ($1+3$)-dimensions, we define
\[A_R := \{ x\in\rt:R<|x|<2R\} \ \ (R>2), \qquad A_{R=1} := \{ |x|<2 \}.
\] 
Given a subinterval $I$ of $[0,\iy)$,
\begin{equation}
\begin{split}
 \| \phi\|_{LE(I)} &:= \sup_R  \| \la r\ra^{-\frac12} \phi\|_{L^2 (I\times A_R)},\\
  \| \phi\|_{LE^1(I)} &:= \| \de \phi\|_{LE(I)} + \| \la r\ra^{-1} \phi\|_{LE(I)},\\
 \| f\|_{LE^*(I)} &:= \sum_R  \| \la r\ra^{\frac12} f\|_{L^2 (I \times A_R)}.
\end{split} 
\end{equation}
We also define 
\[
\begin{split}
  \| \phi\|_{LE^{1,k}(I)} &= \sum_{|\alpha| \leq k} \| \partial^\x \phi\|_{LE^1(I)} \\
  \| \phi\|_{LE^{0,k}(I)} &= \sum_{|\alpha| \leq k} \| \partial^\x \phi\|_{LE(I)},\\
  \| f\|_{LE^{*,k}(I)} &=  \sum_{|\alpha| \leq k}  \| \partial^\alpha f\|_{LE^{*}(I)}.
\end{split}  
\]
For any norm, an omission of $I$ will denote $I:=[0,\iy)$.

\begin{definition}[Local energy decay] \label{LED}
We say that the solution to \eqref{eq:problem} satisfies the local energy decay
estimate if the following estimate holds in $[0,\iy)\times\rt$:
\begin{equation}
 \| \p\|_{LE^{1,k}} 
\ls_k\|\de \p(0)\|_{H^k} + \|f\|_{LE^{*,k}}, \qquad k \geq 0
\label{eq:LED}\end{equation}
\end{definition}

Let $\chi(x)$ be a compactly supported and smooth function equalling 1 in a neighbourhood of the trapped set. We define a weaker version each of the $LE^1$ norm that excises the trapped set region when evaluating $\de \phi$ in $LE$ norm. We also define the attendant dual weak norm.
\begin{align*}
\|\p\|_{LE^1_w(I)} &:= \|(1-\chi)\de \p\|_{LE(I)}+\| \jr\inv \phi\|_{LE(I)}, \ \ \ \|\p\|_{LE^{1,k}_w(I)} := \sum_{|\x|\le k} \|\pa^\x \p\|_{LE^1_w(I)} \\
\|f\|_{LE^*_w(I)} &:= \|f\|_{LE^*(I)} + \|\chi \de f\|_{\lt(I)L^2}, \hspace{1.75cm} \|f\|_{LE^{*,k}_w(I)} := \sum_{|\x|\le k} \|\pa^\x f\|_{LE^{*}_w(I)} %
\end{align*}

We assume that \eqref{eq:problem} satisfies the following weak version of local energy decay \cref{definition:weak}, as expressed by the following bounds: 
\begin{definition}[Weak local energy decay] \label{definition:weak}
	 		We say \eqref{eq:problem} satisfies the weak local energy decay estimate if for any real $T_0\ge0$ and any integer $k\ge0$ 
				\begin{equation} \label{eq:WLED}
	 				\| \p\|_{LE^{1,k}_w[T_0,\iy)} \ls_k \|\de\p(T_0)\|_{H^{k}} + \|f\|_{LE^{*,k}_w[T_0,\iy)}.
				\end{equation}
\end{definition}

\begin{remark}[Loss of two derivatives in the inhomogeneity] %
Combining the $k$ and $k+1$ cases of \cref{eq:WLED} implies
\begin{equation}\label{wled implication}
\|\phi\|_{LE^{1,k}[T_0,\iy)} \ls_k \|\de \p(T_0)\|_{H^{k+1}} + \|f\|_{LE^{*,k+2}[T_0,\iy)}.
\end{equation}
Notice that the right-hand side must have $k+2$ derivatives falling on $f$, since the weak dual norm loses one derivative (at least on $\supp\chi$), and we have applied the $k+1$ case. 
\end{remark}

\begin{definition}[Commuting vector fields and function classes $S^Z$]\label{def:cvf}
In $\R^{1+3}$, we consider the three (ordered) sets 
\[
\partial := (\partial_t, \partial_1,\pa_2,\pa_3), \qquad \Omega := (x^i \partial_j -
x^j \partial_i), \qquad S := t \partial_t + \sum_{i=1}^3x^i \partial_{i},
\]
which are, respectively, the generators of translations, rotations and scaling. We set
\[Z := (\pa,\Omega,S) \] and we define the function class $$S^Z(f)$$ to be the collection of functions $g : \R \times \rt \to \R$ such that
\[ |Z^J g(t,x)| \ls_J |f| \]
whenever $J$ is a multiindex. We will frequently use $f = \jr^k$ for some real $k$.
\end{definition}

\begin{remark}[Instances in which weak local energy decay holds]
Weak local energy decay is known to hold in the Schwarzschild space-time and the Kerr space-time with small $0 \leq |a| \ll M$, where the parameter $M$ denotes the mass of the black hole and the parameter $a$ denotes the angular momentum per unit mass (thus $aM$ denotes the angular momentum of the black hole); more can be found in \cref{ss:Local energy decay estimates}. %

Examples where the assumptions we make on $\p$ in this paper are actually satisfied include the following situations: 
\begin{itemize}
\item
the case with small (meaning $O(\eps)$ in all compact regions for the function and all its derivatives) and asymptotically flat perturbations $h \in S^Z(\eps \jr^{-1-\xs})$ and a small potential $V\in S^Z(\eps\jr^{-2-\xd})$ for arbitrary real numbers $\xd, \xs >0$.  See \cite{L}. 
\item
The situation analysed in \cite{MST}, which proves local energy decay estimates for solutions to scalar wave equations on nontrapping, asymptotically flat space-times (in particular large perturbations of Minkowski space-time).
\end{itemize}
\end{remark}

\begin{remark}[Relation between weak local energy decay and stationary local energy decay]
The problem \cref{main} would be said to satisfy stationary local energy decay estimates (for derivatives) if for any interval $[T_1,T_2]$ and any integer $k\ge0$, we have
\[
\|\p\|_{LE^{1,k}[T_1,T_2]}\ls_k\sum_{j=1}^2\|\de \p(T_j)\|_{H^k} + \|f\|_{LE^{*,k}_\text{w}[T_1,T_2]} + \|\pat \p\|_{LE^{0,k}[T_1,T_2]}.	
\] Notice that we allow for $f$ in the weak dual norm $LE^*_\text{w}$, rather than in the usual dual norm $LE^*$. This differs from the definition of stationary local energy decay in \cite{MTT}.

There is an analogous version of the stationary local energy decay estimates for commuting vector fields (see \cref{def:cvf}; they are also simply called \textit{vector fields}) of $\phi$.

One can prove the following: if the weak local energy decay \cref{definition:weak} for derivatives of $\phi$ holds and $\pat$ is timelike on any trapped region that may exist, then $\phi$ in fact also satisfies stationary local energy decay estimates for \textit{vector fields of} $\phi$. In this sense, the weak local energy decay is a weaker assumption than stationary local energy decay.

In this paper, we do not make the assumption of stationary local energy decay, and we do not assume that $\pat$ is timelike on any trapped region that may exist. This is in contrast to the assumptions made in \cite{MTT}, where the authors do make these assumptions. Our argument in this paper generalises the argument in \cite{MTT} in this regard, and also in the decay rates of our coefficients, the kinds of coefficients considered in the wave operator $P$, and the support of the initial data. We consider solutions with non-compactly-supported initial data, in fact solutions with initial data in a weighted $\lt$ space (see \cref{thm:main} for the main theorem and assumptions).
\end{remark}

\begin{definition}
We define $S^Z_\text{cone}(f)$ to be the collection of $g$ such that $|Z^J g| \ls |f|$ in 
$\{ t/2 \leq r \leq 3t/2\}$. 
	Thus $S^Z(f)\subsetneq S^Z_\text{cone}(f)$. 
 We define $S^Z_\text{int}(f)$ to be the collection of $g$ such that $|Z^J g| \ls |f|$ in $\{ r < t/2 \}$. 
	We define $S^Z_\text{radial}(f) := \{ g \in S^Z(f) : g \text{ is spherically symmetric} \}$. 

Let $\| \cdot\|$ be any norm used in this paper. Given any nonnegative integer $N\ge0$, we write $\|g_{\le N}\|$ to denote $\sum_{|J|\le N} \|g_J\|$. (See also \cref{norm def}.)
\end{definition}   

Our main result, \cref{thm:main}, is a pointwise decay estimate for the solution to the following equation: 
\begin{equation}   \label{main}
\begin{cases}
  P\phi(t,x) =  0 \qquad (t,x) \in (0,\infty) \times \R^3, \  P \text{ given in }\eqref{P def}  \\
(\p(0,x), \bar N\p(0,x)) = (\p_0(x),\p_1(x)) 
\end{cases}
\end{equation}
where $g$ is a Lorentzian metric and $g^\xo,A,B,V$ are functions satisfying the following conditions:

\begin{theorem}[Main theorem] \label{thm:main}
Let $m\ge0$ be an integer and let $N$ be a sufficiently large integer relative to $m$, or $N\gg m$.
Let $\gab(t,x)$ be a Lorentzian metric such that for all $t_0\ge0$ the level sets $\{ t =t_0\}$ are space-like, 
and let $h := g-m$ with $m$ denoting the Minkowski metric. Let $u := t-r$ and $v:=t+r$. 
	Assume that $\p$ solving \eqref{main} satisfies 
the weak local energy decay \eqref{eq:WLED}, and that $\phi_0 \in \lt(\R^3)$. 

\begin{enumerate}
\item
Suppose that for some real $0<\xs, \xd, \xd'<\iy$,
$$h \in S^Z(\jr^{-1-\xs})$$ 
$$A\in S^Z(\jr^{-1-\xs})$$
$$\pat A\in S^Z\lr{\jv\ju\inv\jr\inv \jr^{-1-\xs} } \cap S^Z_\text{cone}(\jr^{-1-\xs})$$
$$\pa A\in S^Z_{int}(\jr^{-2-}) \cap S^Z_{cone}(\jr^{-2-})$$
$$B\in S^Z(\jr^{-1-\xs})$$
$$\pat B\in S^Z(\jr^{-2-\xs})$$
$$V \in S^Z(\jr^{-2-\xd})$$
$$g^\xo \in S^Z_\text{radial}(\jr^{-2-\xd'})$$

Then %
\begin{equation}\label{eq:main bound 1}
|\pm(t,x)| \ls \f1{\la v\ra\la u\ra^{1+\min(\xs,\xd,\xd')}}\|\jr^{\f12+\min(\xs,\xd,\xd')} \de\p_{\le N}(0)\|_{\lt(\rt)}
\end{equation}
\begin{equation}
|\pa\pm(t,x)| \ls \f1{\jr \ju^{2+\min(\xs,\xd,\xd')}}\|\jr^{\f12+\min(\xs,\xd,\xd')} \de\p_{\le N}(0)\|_{\lt(\rt)}
\end{equation}
\begin{equation} %
|\pat\pm(t,x)| \ls \f1{\jv \ju^{2+\min(\xs,\xd,\xd')}}\|\jr^{\f12+\min(\xs,\xd,\xd')} \de\p_{\le N}(0)\|_{\lt(\rt)}
\end{equation}
\begin{equation}
|\pa^2\pm(t,x)| \ls \f\jv{\jr^2 \ju^{3+\min(\xs,\xd,\xd')}}\|\jr^{\f12+\min(\xs,\xd,\xd')} \de\p_{\le N}(0)\|_{\lt(\rt)}
\end{equation}
\item
If in addition to the assumptions in part (1),
	$$\pat h\in S^Z\lr{ \jr^{-2-\xs} }$$
$$\pat^2 h\in S^Z_\text{cone}\lr{ \ju^{-2} \jr^{-1-\xs} }$$
$$A\in S^Z(\jr^{-2-\xs})$$ 
$$\pat A\in S^Z\lr{ \jv\jr\inv\ju\inv \jr^{-2-\xs} }$$ 
$$B\in S^Z(\jr^{-2-\xs} )$$
$$\pat B\in S^Z(\jr^{-3-\xs})$$
then the solution to \eqref{main} satisfies
\begin{equation}\label{eq:main bound 2}
|\pm(t,x)| \ls \f1{\la v\ra\la u\ra^{1+\min(1+\xs,\xd,\xd')}}\|\jr^{\f32+\min(\xs,\xd,\xd')} \de\p_{\le N}(0)\|_{\lt(\rt)}
\end{equation}
\begin{equation}
|\pa\pm(t,x)| \ls \f1{\jr\ju^{2+\min(1+\xs,\xd,\xd')}}\|\jr^{\f32+\min(\xs,\xd,\xd')} \de\p_{\le N}(0)\|_{\lt(\rt)}
\end{equation}
\begin{equation}
|\pat\pm(t,x)| \ls \f1{\jv\ju^{2+\min(1+\xs,\xd,\xd')}}\|\jr^{\f32+\min(\xs,\xd,\xd')} \de\p_{\le N}(0)\|_{\lt(\rt)}
\end{equation}
\begin{equation}
|\pa^2\pm(t,x)| \ls \f\jv{\jr^2\ju^{3+\min(1+\xs,\xd,\xd')}}\|\jr^{\f32+\min(\xs,\xd,\xd')} \de\p_{\le N}(0)\|_{\lt(\rt)}
\end{equation}
\end{enumerate}
\end{theorem}

\begin{remark}We make some remarks supplementing the main theorem.
\begin{itemize}
\item
Second order operators that have spherically symmetric coefficients of the form $1/r$, at least away from the origin, 
are covered, i.e. included, by the definition of our operator $P$ in \cref{P def}.\footnote{Indeed, \cref{P def} covers coefficients that have the following form away from the origin: $1/r^{a}$, $a \in \R_{>0}$.} This appears in some equations of physical interest, such as in general relativity.
\item
The argument shown in this paper straightforwardly yields a proof of a more general version of \cref{thm:main} which assumes more general decay rates on $A$ and $B$. 
	Namely, given any real $\xs',\xs''>0$, for part (1) of \cref{thm:main} (and similarly for part (2)), if
$$A\in S^Z(\jr^{-1-\xs'})$$
$$\pat A\in S^Z\lr{\jv\ju\inv\jr\inv \jr^{-1-\xs'} } \cap S^Z_\text{cone}(\jr^{-1-\xs'})$$
$$B\in S^Z(\jr^{-1-\xs''})$$
$$\pat B\in S^Z(\jr^{-2-\xs''})$$
(in addition to the assumptions on $h,g^\xo$ and $V$ in part (1), as well as the assumption on the generic derivative $\pa A$)
then the same arguments in this paper automatically give, for instance, with part (1) assumptions,
$$|\pm(t,x)| \ls \f1{\la v\ra\la u\ra^{1+\min(\xs,\xd,\xd',\xs',\xs'')}}\|\jr^{\f12+\min(\xs,\xd,\xd',\xs',\xs'')} \de\p_{\le N}(0)\|_{\lt(\rt)},$$
and the corresponding bounds also hold for $\pat \pm, \pa\pm$, and so on. %

For simplicity of presentation, in this paper we restrict to the case $\xs = \xs' = \xs''$.
\item
In item (2) of \cref{thm:main}, one class of examples of metrics $\gab$ satisfying the conditions given are the stationary metrics $g$, that is, those with stationary component $$h =h(x).$$ By substituting the natural number values $\xd \ge 1, \xd\in\N$ and $\xs \ge2, \xs\in\N$, this special case of item (2) of \cref{thm:main} recovers a similar result as the main theorem in \cite{Mor}.
\item
If $\hab \in S^Z(\jr^{-q})$ for some $q>0$, then $\sqrt{|g|}\hab \in S^Z(\jr^{-q})$. This is a consequence of the product rule and the assumption that $-q < 0$. 
Thus \cref{thm:main} also holds if $\paa\gab\pab$ is replaced by the geometric wave operator 
\[
\Box_g=\frac{1}{\sqrt{|g|}} \paa \sqrt{|g|} \gab \pab, \quad |g|:=|\det \gab|.
\] 
\item One statement this theorem says is that if local energy estimates even with derivative loss is assumed (see \cref{wled implication}), then one can obtain the pointwise bounds in \cref{thm:main}. 
\item
For a first reading, since $ S^Z\lr{\jv\ju\inv\jr\inv \jr^{-1-\xs} } \cap S^Z_\text{cone}(\jr^{-1-\xs})\subset S^Z(\jr^{-2-\xs})$, the reader may wish to keep in mind that $\pat A\in S^Z(\jr^{-2-\xs})$ for part 1 of \cref{thm:main}. 
\end{itemize}
\end{remark}

\begin{remark}[Black hole spacetimes]
All the arguments in this paper can be adapted to the exterior of a ball and hence the proofs in this paper can be applied in the case of black hole spacetimes.
\end{remark}

\subsection{Local energy decay estimates}\label{ss:Local energy decay estimates} The first instance of a local energy estimate was obtained by Morawetz for the Klein-Gordon equation in \cite{M}. Some other work on local energy decay estimates and their applications can be found in, for instance, \cite{Strauss,Al,KSS,KPV,MS,MT,SmSo,St}. For local energy decay estimates for small and time dependent long range
perturbations of the Minkowski space-time, see for instance \cite{Al}, \cite{MT2}, \cite{MS} for time dependent
perturbations, as well as, e.g., \cite{B}, \cite{BH}, \cite{SW} for time independent,
nontrapping perturbations. There is a related family of local energy decay estimates for the Schr\"odinger equation as well.

For Schwarzschild metrics, trapping at the event horizon was shown to be trivial due to an effect guaranteeing energy decay along the trapped rays called the redshift effect. On the other hand, for Kerr metrics, a local energy estimate with derivative loss on the trapped set is often introduced. \cref{definition:weak} includes this loss.

	For large perturbations of the Minkowski metric, if one assumes the absence of trapping then local energy estimates can still hold; see for instance \cite{BH, MST}. 
	For weak enough trapping, \cref{definition:weak} has been established; see for instance \cite{BCMP,Chr,NZ,WZ}. 
	If one assumes absence of trapping, then %
	\cref{LED} holds;
		with trapping \cref{LED} cannot hold, see \cite{Ral, Sb}. 
	With sufficiently strong trapping, even \cref{definition:weak} fails, see \cite{DMST}.

Weak local energy decay for the Schwarzschild metric was established in \cite{B, DR2, MMTT}. For the Kerr metric with low angular momenta, weak local energy decay estimates were proved in \cite{B, DR2, DRSR}. 

The local energy estimate for Kerr spacetimes with small angular momenta was proven in \cite{TT} (see also \cite{AB} and \cite{DaRoNotes} for related work), for large angular momentum $|a|<M$ in \cite{DRSR}, and for extremal Kerr $|a|=M$ in \cite{Ar}.

\subsection{Pointwise decay and asymptotic behaviour}
It is well-understood that local energy decay in a compact region on an asymptotically flat region implies pointwise decay rates that are related to how rapidly the metric coefficients decay to the Minkowski metric; see, for example, the works \cite{Tat, MTT,OS,Mos,AAG1,AAG2,Mor, MW, Hin2, LiT2,L,Loo22}. 
Similar results as this paper, for stationary spacetimes, were shown in \cite{MW} concurrently using spectral theory techniques, with particular focus on analyzing the resolvent at low frequencies.

Local energy decay is also involved in proving scattering, another type of asymptotic behaviour, on variable-coefficient backgrounds. In particular, they imply Strichartz estimates on certain variable-coefficient backgrounds, see \cite{MT}. The article \cite{LT} used local energy decay to prove scattering for the version of the problem \eqref{main} without the potential $V$ and first-order terms $A$ and $B$, although the argument extends straightforwardly to the problem including $V,A$ and $B$ defined above. 

In the case of the Schwarzschild metric, Price \cite{Pri} conjectured that the solution to the wave equation decays at the rate $t^{-3}$ within any compact region; this rate was shown to hold for a variety of spacetimes, including Schwarzschild and Kerr spacetimes with small angular momenta---see \cite{DSS,Tat,MTT}. 

\subsection{The main ideas of the proof}  Aside from the standard tools of Sobolev embedding, albeit exploited primarily in dyadic conical subregions, when proving pointwise bounds we take advantage of the reduction to $1+1$ dimensions in spherical symmetry---called the ``one-dimensional reduction''---and the positivity of the fundamental solution to the $1+3$ dimensional wave equation. This not only provides a simple setting for the analysis but also allows us to ``absorb'' pointwise decay from the vector fields of the coefficients $h, V,$ and so forth, and transfer them to the decay of the solution $\p$ or its vector fields. 
	In this way, gradual improvements, starting from an initial decay estimate \cref{u/v decay}---obtained from only Sobolev embedding and integrated local energy decay---are possible, with the improvements arising from the positivity of $\xs, \xd$ and $\xd'$. 

A little more precisely, for components of the wave equation that contain a derivative structure, we analyse them separately: in a neighbourhood of the light cone $\{r=t\}$ (see \cref{sec:cone bounds}); and in all other regions. However, for components of the wave equation that do not contain a derivative structure, we need not make this distinction. 

\subsection{Summary of sections} 
 
In \cref{sec:notation}, we define notation and conventions for the rest of the paper.

In \cref{sec:commuting}, we commute $P$ with vector fields and prove (weak) local energy estimates for vector fields.

In \cref{sec:initialests}, we prove Sobolev embedding estimates and obtain an initial pointwise decay estimate.       We connect pointwise bounds to $L^2$ estimates and norms, thereby connecting local energy decay to pointwise bounds.

In \cref{sec:der ests}, we prove that derivatives of vector fields of the solution decay better at the cost of applying more vector fields.

In \cref{sec:setup}, we define more notation that will be used for the pointwise decay iteration, which occupies the remainder of the paper. We also prove certain lemmas used in the iteration.

In \cref{sec:exte}, we prove the upper bound in $\{r>t+1\}$ for components of the solution away from the cone.

In \cref{sec:r to t}, we show how to convert a decay rate of $\jr^{-p}$ for the solution $\p$ and its vector fields to $\la t+r\ra^{-p}$ for $p\le1$.

In \cref{sec:inte}, we prove the upper bound in $\{r<t\}$ for components of the solution away from the cone.

In \cref{sec:cone bounds}, we prove the upper bound for components of the solution near the cone.

\section{Notation and conventions}\label{sec:notation}

\subsection{Notation for dyadic numbers and conical subregions}
We work only with dyadic numbers that are at least 1. We denote dyadic numbers by capital letters for that variable; for instance, dyadic numbers that form the ranges for radial (resp. temporal and distance from the cone $\{|x|=t\}$) variables will be denoted by $R$ (resp. $T$ and $U$); thus $$R,T, U\ge 1.$$ 
	We choose dyadic integers for $T$ and a power $a$ for $R,U$---thus $R = a^k$ for $k\ge1$--- different from 2 but not much larger than 2, for instance in the interval $(2,5]$, such that for every $j\in\N$, there exists $j'\in\N$ with 
\begin{equation}\label{dyadic numbering}
a^{j'} = \f38 2^j.
\end{equation} 
 
\subsubsection{Dyadic decomposition}
We decompose the region $\{r\le t\}$ based on either distance from the cone $\{r=t\}$ or distance from the origin $\{r=0\}$. We fix a dyadic number $T$. 
\begin{align*}
C_T &:= \begin{cases}
\{ (t,x) \in [0,\iy) \times \rt : T \leq t \leq 2T, \ \ r \leq t\} & T>1 \\
\{ (t,x) \in [0,\iy) \times \rt : 0 < t < 2, \ \ r \leq t\} & T=1
\end{cases} \\
C^R_T &:=\begin{cases}
C_T\cap \{R<r<2R\} & R>1\\
C_T\cap \{0 < r < 2\} & R=1
\end{cases}\\
C^U_T &:=\begin{cases}
\{ (t,x) \in [0,\iy) \times \rt  :  T\le t\le 2T\} \cap \{U<|t-r|<2U\} & U>1\\
\{ (t,x) \in [0,\iy) \times \rt  :  T\le t\le 2T\} \cap \{0< |t-r|<2\} & U=1
\end{cases} %
\end{align*}
If a need arises to distinguish between the $R=1$ and $U=1$ cases, we shall write $C^{R=1}_T$ and $C^{U=1}_T$ respectively. Note that $|C^R_T| \sim (R^3T)^{1/2}$ and $|C^U_T| \sim (T^3 U)^{1/2}$. 
	We define
\begin{equation*}
C_T^{<3T/4} := \bigcup_{R < 3T/8} C_T^R.
\end{equation*}

Now letting $R > T$, we define
\begin{align*}
C^T_R &:= \{ (t,x) \in [0,\iy) \times \rt  :  r \ge t, T \le t\le 2T, R \le r\le 2R, R \le |r-t| \le 2R\}
\end{align*} 
Note that $|C^T_R| \sim R^2$, as can be seen in the $|r-t|$ and $r$ directions. 

	$C_T^R, C_T^U$ and $C^T_R$ are where we shall apply Sobolev embedding, which allows us to obtain pointwise bounds from $L^2$ bounds. 

\subsubsection{Enlargements of sets}\label{enlarge}
Given any subset of these conical regions, a tilde atop the symbol $C$ will denote a slight enlargement of that subset; for example, $\crtt$ denotes a slightly larger set containing $\crt$. 

\subsection{More notation for vector fields}
Beyond \cref{def:cvf}, we now define more notation for vector fields. 

\subsubsection{Subscripts on functions will denote vector fields.}
  Given a nonnegative integer $m$ and a triplet $J = (i,j,k)$ of multi-indices $i$, $j$ and $k$ for $(\pa,\xO,S)$---by this we mean $\pa^i\xO^jS^k$---we
denote $|J| = |i|+4|j|+10k$. See \cref{vf defn}. 

The coefficient 10 in front of $k$ arises because of the fact $[P,S] - 2P - s_{2+\xd'} \xO^2 \in \mathcal C$, where $\mathcal C$ is the class of operators defined in \cref{commdef}. In particular the presence of $\xO^2$ as well as loss of derivative considerations (a price of losing two derivatives if one wants to control the full $LE^1$ norm---see \cref{wled implication}) for the inhomogeneity in the weak local energy decay \cref{definition:weak} leads to the count $10 = 2 + 2\cdot4$. 
If $g^\xo=0$, then we would count each $S$ in the same way we would count $\pa^2$, i.e. two derivatives. \ We put in place these differences in these numerical weights for $i,j,$ and $k$ (respectively: 1, 4, and 10) because of the trapped set.

We denote
\begin{equation}\label{vf defn}
\p_{J} := Z^J\p := \partial^i \Omega^j S^k \p, 
\end{equation}
\[ 
\p_{\leq m} := (\p_{J})_{|J|\le m}, \ \ \  \phi_{m_1 \leq \cdot \le m_2} := (\p_{J})_{m_1 \leq |J|\le m_2}, \ \ \ \p_{=m} := (\p_J)_{|J|=m}
\]
\[
\pa^{\le m} \phi := (\pa^{i'}  \phi)_{|i'| \leq m}, \ \ \ \pa^{= m} \phi := (\pa^{i'}  \phi)_{|i'| = m}%
\]
Furthermore, by $Z^{=m}\phi$ we mean $\p_{=m}$, and so on. 
	We write $J_1 \leq J_2$ to mean
\[
i_1 \leq i_2, \qquad j_1\leq j_2, \qquad k_1\leq k_2,
\]
and $J_1 < J_2$ if at least one of the inequalities above is strict. 
	If $I$ is a multiindex of order $\ell$ and $n$ an integer, by $I+n$ we mean 
	$$\{ I + J : |J| = n, J \text{ is an }\ell \text{-multiindex} \}.$$ 
	Given a multiindex $K$, we define 
\[
\p_{\le K} := (\p_J)_{J \le K}.
\]

\subsection{Notation for the symbols $n$ and $N$} \label{subsec:N}
Throughout the paper the integer $N$ will denote a fixed and sufficiently large positive number%
signifying the highest total number of vector fields that will ever be applied to the solution $\p$ to \eqref{main} in the paper. 

We use the convention that the value of $n$ may vary by line.%

\subsection{The use of the tilde symbol} If $\xS$ is a set, we shall use $\ti\xS$ to indicate a slight enlargement of $\xS$, and we only perform a finite number of slight enlargements in this paper to dyadic subregions. The symbol $\ti\xS$ may vary by line. 

If $f$ is a function, we shall typically use $\ti f$ to denote commuting vector fields applied to $f$.

\subsection{Notation for implicit constants} \label{subsec:implicit constants}
We write $X\ls Y$ to denote $|X| \leq CY$ for an implicit constant $C$ which may vary by line. Similarly, $X \ll Y$ will denote $|X| \le c Y$ for a sufficiently small constant $c>0$. In this paper, all implicit constants are allowed to depend on the dimension and the initial data $\phi_{\leq N}[0]$, for a fixed $N\in \N$ that is sufficiently large.

\begin{definition}\label{norm def}
Given a real number $t \geq 0$, a norm or square of a norm $\calF[\cdot](t)$ (including the absolute value), a function $f$, a multiindex $J$, and an integer $k$, we let 
\begin{align*}
\calF[f_{\le J}](t) &:= \sum_{\text{multiindices }I\le J} \calF[f_I](t)\\
\calF[f_{\le k}](t) &:=\sum_{\text{multiindices }I : |I|\le k}\calF[f_I](t)
\end{align*}
\end{definition}

\subsection{Other notation}

If $t\geq0$ is a real number, let $$\xS_t:=\{ (t,x) : x\in\rt\}$$ denote the constant time $t$ slice. 

If $x =(x^1,x^2,x^3)\in\R^3$, we write 
\begin{align*}
r &:= \lr{ \sum_{i=1}^3 (x^i)^2 }^{1/2}, \\
u&:= t-r,   \\   %
v&:= t+r.  %
\end{align*}
We write $\bo := -\pat^2+\Delta$. %

We write $$s_{q}$$ to denote element of $S^Z(\jr^{-q})$. $q$ will denote a nonnegative number. %

\section{Commuting with vector fields, and weak local energy decay for vector fields}\label{sec:commuting}

\begin{remark}
Let $w$ be a sufficiently smooth function. Then
\begin{equation}\label{D decomp}
\pa w\in S^Z(\jr\inv)  \bar Zw + \mu S^Z(1) |\pat w| \text{ if }r \ge t/2
\end{equation}
with $\mu=0, \bar Z=\xO$ for angular derivatives $\pa_\xo w$ on the left-hand side, and $\mu=1, \bar Z=S$ for the radial derivative $\pa_r w$ on the left-hand side. 
\end{remark}

\

We define $\mathcal C$ to be the collection of real linear combinations of the operators
\begin{equation}\label{commdef}
 \pa s_{1+q'} \pa, \  s_{1+q'} \pa \pa, \ s_{2+q'} , \ \pa s_{1+q'},  \ s_{1+q'} \pa
\end{equation} 
where $q'>0$ is a number which depends on the assumptions made about the coefficients $h,g^\xo,V,A,$ and $B$ in \cref{thm:main}. That is, schematically, $\mathcal C = \{  \pa s_{1+q'} \pa + s_{1+q'} \pa \pa + s_{2+q'} + \pa s_{1+q'}+ s_{1+q'} \pa \}$.

\begin{lemma} \label{lem:comm0}
Let $w$ be a sufficiently smooth function. Given %
$J$ and $k\ge0$, there are some operators $\dot C\in \mathcal{C}$ such that
\begin{equation}\label{rot and scal comm}
\xO^J (S+2)^k Pw = P \xO^J S^k w + \dot C w_{\le 4(|J|-1) + 10k}
\end{equation}
where we adopt the following conventions:  we interpret $\dot C w_{\le 4(|J|-1) + 10k}$ as a sum, and subscripts with negative real value denote the zero multiindex.
\end{lemma}

\noindent \textit{Proof (sketch).} 
By the assumptions in the main theorem,
\begin{equation}\label{commfact1}
[P,\pa]  \in \mathcal{C}.
\end{equation}
\begin{equation}\label{commfact2}
[P,\xO]  \in \mathcal{C}.
\end{equation}
\begin{equation}\label{commfact3}
[P,S] - 2P - s_{2+\xd'} \xO^2  \in \mathcal{C}.
\end{equation}
		One uses \cref{commfact1,commfact2,commfact3} and proves the result by mathematical induction. We omit the details of the proof, except for the following observation.
	Starting from $\xO^J(S+2)^k P$ and then commuting the vector fields with $P$, then other than $P\xO^JS^k$, the terms with the highest vector field count (assuming $g^\xo$ is not the zero function) are those of the form  
	$$\dot C \bar Z^{= |J| + k-1} w, \ \ \bar Z\in\{\xO,S\}, \ \ \dot C\in \mathcal C;$$
	 more specifically, those of the form $\dot C \xO^{|J|-1}S^k$. This explains the subscript $4(|J|-1)+10k$. 

\begin{lemma}%
\label{lem:comm}
Given the assumptions in either part 1 or part 2 of \cref{thm:main}, there exists a positive real number $q'>0$ such that for any multiindex $J$,
$$|P\p_J| \ls \f{|\p_{\leq |J|-1}|}{\jr^{2+q'}} + \f{|\de\p_{\leq |J|}|}{\jr^{1+q'}} + |(P\p)_{\le |J|}|.$$
\end{lemma}

\begin{proof}
There is a constant $q'>0$ such that the operator $P$ can be written schematically as
$P=\Box+\pa s_{1+q'} \pa + s_{1+q'} \pao^2 + s_{2+q'} + s_{1+q'}\pa + \pa s_{1+q'}$. We have $[Z,\pa] = c\pa$ schematically, for some real number $c$ depending on $Z$. 

For terms of the form $(\pa\ti A) \ti \p$, where $\ti A, \ti\p$ denote possible vector fields %
of $A,\p$, we apply the assumption 
$$\pa A\in S^Z_{int}(\jr^{-2-}) \cap S^Z_{cone}(\jr^{-2-})$$
on generic derivatives $\pa A$ from part 1 of \cref{thm:main} in $\{ r<3t/2\}$, and the assumption on $\pat A$ and \eqref{D decomp} in $\{ r \ge 3t/2\}$, giving a contribution of the form $\jr^{-2-q'}|\p_{<|J|}|$.
	For part 2, on the other hand, we in fact need not look at $r<3t/2$ and $r\ge 3t/2$ separately, because the statement $\pa A\in S^Z(\jr^{-2-})$ is already trivially satisfied for any $(t,r)$-pair given the assumption on $A$. 
	
	We include the terms arising from $g^\xo\Delta_\xo$ together with the $\jr^{-1-}|\de \p_{\le|J|}|$ term.
	The rest is clear, and the claim follows.
\end{proof}

We recall the weak local energy decay estimate $\| \p\|_{LE^{1,k}_w} \ls_k \|\de \p(T_0)\|_{H^{k}} + \|f\|_{LE^{*,k}_w},$
 which can be rephrased as 
\begin{align*}
\sum_{|\x|= k+1} \|(1-\chi)\pa^\x\p\|_{LE[T_0,\iy)} 
	&+ \|\jr\inv \p\|_{LE[T_0,\iy)} 
	+ \sum_{1 \le|\xg|\le k} \|\pa^\xg\p\|_{LE[T_0,\iy)} \\
	&\ls_{k,\chi} \|\de\p(T_0)\|_{H^{k}} + \|f\|_{LE^{*,k}_w[T_0,\iy)}.
\end{align*}

\begin{proposition}[Weak local energy decay for vector fields] \label{wled for vf}
Let $\p$ be any smooth-enough function solving \cref{main} and satisfying \cref{definition:weak}. 
	Then for any natural number $m\geq 0$, 
\begin{align}\label{prop.claim}
\begin{split}
\|\phi_{\le m}\|_{LE^1} \ls \|\de \pmo(0)\|_\lt + \|f_{\le m+2}\|_{LE^*}.
\end{split}
\end{align}
\end{proposition}

\begin{proof} We prove \eqref{prop.claim} by induction. 

The base case
\begin{align*}
\|\p\|_{LE^1}\ls \|\de\p_{\le1}(0)\|_\lt + \|f_{\le2}\|_{LE^*}
\end{align*}
is simply given by combining \cref{definition:weak} at $k=0$ and $k=1$, which yields
$$\|\p\|_{LE^1}\ls \|\de\p(0)\|_{H^1} + \|\pa^{\le1}f\|_{LE^*} + \|\chi \pa^{=2} f\|_{\ltlt},$$
which is clearly bounded by 
$$\|\de\p(0)\|_{H^1} + \|\pa^{\le2}f\|_{LE^*} \le \|\de\p_{\le1}(0)\|_{\lt} + \|f_{\le2}\|_{LE^*}.$$

Next, we use \cref{lem:comm0}. Let $|(I,J,k)| = m$. 
\begin{align*}
\|\p_{(I,J,k)}\|_{LE^1} &\ls \|\de \xO^JS^k\p(0)\|_{H^{|I|+1}} + \|\xO^JS^kf\|_{LE^{*,|I|+2}} + \|[P,\xO^JS^k]\p\|_{LE^{*,|I|+2}}\\
&\ls \|\de\p_{\le m+1}(0)\|_\lt + \|f_{\le m+2}\|_{LE^*} + \|[P,\xO^JS^k]\p\|_{LE^{*,|I|+2}}\\
&\ls \|\de\p_{\le m+1}(0)\|_\lt + \|f_{\le m+2}\|_{LE^*} + \|\jr^{-1-}\de \p_{\le m-2}\|_{LE^*} + \|\jr^{-2-} \p_{\le m-2}\|_{LE^*} \\
&\ls \|\de\p_{\le m+1}(0)\|_\lt + \|f_{\le m+2}\|_{LE^*} + \| \p_{\le m-2}\|_{LE^1}\\
&\ls \|\de\p_{\le m+1}(0)\|_\lt + \|f_{\le m+2}\|_{LE^*}
\end{align*} In transitioning from the second line to the third line, we used \cref{rot and scal comm}. 
The third line follows by the assumption that $\xO$ counts for four partial derivatives.%
 The final line follows by the induction hypothesis.
\end{proof}

\begin{remark}
The above proof extends to time intervals $[T_1,\iy)$, $T_1\geq 0$. (The proof above assumes $T_1=0$.) %
The estimate is
$$\|\phi_{\le m}\|_{LE^1[T_1,\iy)} \ls \|\de \pmo(T_1)\|_\lt + \|f_{\le m+2}\|_{LE^*[T_1,\iy)}.$$
\end{remark}

\section{Initial $L^\iy$ estimates}\label{sec:initialests}

We now state the Sobolev embedding estimates localised to our selected conical regions. 
\begin{lemma} \label{SobEmbExt}
Let $w\in C^4$. 
\begin{itemize}
\item
For all $T \geq 1$ and $1 \leq U \leq 3T/8$, %
we have  \begin{equation}\label{U-Sob}
 \| w\|_{L^\infty(\cut)} \ls \sum_{i\le 1,j\le 2} \f1{(T^3U)^{1/2}}  \|S^i \Omega^j w\|_{L^2(\ti C^U_T)} +  \lr{\f{U}{T^3}}^\f12 \|\pa_r S^i \Omega^j w\|_{L^2(\ti C^U_T)}.
\end{equation}
\item
For all $T\geq 1$ and $R>T$, we have 
\begin{equation}\label{R-Sob}
  \| w\|_{L^\infty(C^T_R)} 
  \ls \sum_{i\le 1,j\le 2}
\f1{(R^3T)^{1/2}}    \|S^i \Omega^j w\|_{L^2(\ti C^T_R)} +  
\f1{(RT)^{1/2}}   \|\pat S^i \xO^j  w\|_{L^2(\ti C^T_R)}.   %
\end{equation}
\item
For all $T \geq 1$ and $1 \leq R \leq 3T/8$, we have
 \begin{equation}
  \!  \| w\|_{L^\infty(\crt)} 
  \ls \sum_{i\le 1,j\le 2}
\f1{(R^3T)^{1/2}} \|S^i \xO^j w\|_{L^2(\ti C^R_T)} 
   + 
\f1{(RT)^{1/2}} \|\pa_r S^i \xO^j  w\|_{L^2(\ti C^R_T)}.
    \label{in-R-Sob}
\end{equation}
\end{itemize}
\end{lemma}

\begin{proof}

In $\cut$ we make the change of coordinates $t = e^s$ and $|r-t| = e^{s+\rho}$. With this change of coordinates, we are now dealing with a region of size 1 in spherical coordinates including $s$. We have 
$\pa_s = t\pat +r\pa_r =S$ and $\pa_\rho = (r-t)\pa_r$. 
Then we apply the fundamental theorem of calculus in $s$ and also in $\rho$. Finally, we rescale to $\cut$, obtaining \eqref{U-Sob}.

For $C^T_R$, we let $r = e^s$ and $r-t = e^{s+\rho}$. Thus $\pa_s=S$ and $\pa_\rho =(t-r)\pat$. We get 
$$\| w\|_{L^\infty(C^T_R)} 
  \ls \sum_{i\le 1,j\le 2}
\f1{(R^3T)^{1/2}}    \|S^i \Omega^j w\|_{L^2(\ti C^T_R)} +  
    \f{R-T}{(R^3T)^{1/2}}  \|\pat S^i \xO^j  w\|_{L^2(\ti C^T_R)}. $$
 This implies \cref{R-Sob} since $R - T \leq R$. 

For $C^R_T$, we let $t = e^s$ and $r=e^{s+\rho}$. We obtain $\pa_s=S$ and $\pa_\rho=r\pa_r$ and \eqref{in-R-Sob}. 
\end{proof}

\begin{corollary} \label{Sob.emb}
\begin{equation}\label{eq:des}
	\|\phi\|_{L^\iy_{t,x}(\inte)} \ls   \sum_{i \le1,j\leq 2} \f1{T^{1/2}} \| S^i \Omega^j \phi \|_{LE^1_{t,x}(\ti C^{<3T/4}_T)}.
\end{equation}
\end{corollary}

\begin{proof}
By rewriting \eqref{in-R-Sob} in the local energy norm by shifting the $R$ weights around, we obtain \cref{eq:des}.
\end{proof}

\begin{lemma}\label{Hardy}
If $f\in C^1( [0,\iy)_t \times \R^3_x )$, then
\begin{equation}\label{con.Hardy}
\int_{t/2}^{3t/2} \f{f(t,x)^2}{\ju^2} dx \ls \int_{t/4}^{7t/4}|\pa_r f(t,x)|^2 dx + \f1{t^2} \left(\int_{t/4}^{t/2} f(t,x)^2 dx + \int_{3t/2}^{7t/4} f(t,x)^2 dx\right)
\end{equation}
\end{lemma}

\begin{proof}

Let $\chi : [0,\iy) \to [0,1]$ be a cutoff such that $\chi(s) = 1$ for $1/2 \le s \le 3/2$ and 0 when $s \le 1/4$ and $s \ge 7/4$.  We will show that, if $\gamma>-1/2$, and $\gamma \neq 1/2$, then
\begin{align*}
\int\ju^{-2-2\xg} \chi(r/t) f(r,\xo)^2 r^2dr
&\ls  \int\ju^{-2\xg}|\pa_rf(r,\xo)\chi(r/t)|^2 r^2dr \\
	&\quad + \f1{t^2} \int \ju^{-2\xg} |f(r,\xo) \chi'(r/t)|^2 r^2 dr.
\end{align*}
The conclusion follows if we take $\gamma=0$ and integrate over $\omega$.
 
We have
$$f(r,\xo)^2\chi(r/t) - f(7t/4,\xo)^2\chi((7t/4) / t) = -2\int_r^{7t/4} f(\rho,\xo) \chi(\rho/t) \cdot \pa_r ( f(\rho,\xo) \chi(\rho/t) )d\rho.$$
Hence 
$$f(r,\xo)^2 \chi(r/t) r^2 \ls f(7t/4,\xo)^2 \chi(3t/2) t^2 + 2 \int_r^{7t/4} |f(\rho,\xo) \chi(\rho/t) \cdot \pa_r ( f(\rho,\xo) \chi(\rho/t) ) | d\rho$$
Recall that $\chi(7t/4) = 0$. 
We multiply by $\ju^{-2-2\xg}$ and integrate $r$ from $t/4$ to $7t/4$.  This yields
\begin{align*}
\int_{t/4}^{7t/4}\ju^{-2-2\xg} \chi(r/t) f(r,\xo)^2 r^2dr &\ls \int_{t/4}^{7t/4} \ju^{-1-2\xg} |f(r,\xo)\chi(r/t) \pa_r (f(r,\xo)\chi(r/t))| r^2dr
\end{align*}%
By the chain rule, $\pa_r(\chi(r/t)) \ls \chi'(r/t) \cdot \f1t$.
Thus by Cauchy-Schwarz and the chain rule
\begin{align*}
\int_{t/2}^{3t/2}\ju^{-2-2\xg} f(r,\xo)^2 r^2dr
&\ls  \int_{t/4}^{7t/4}\ju^{-2\xg}|\pa_rf(r,\xo)\chi(r/t)|^2 r^2dr \\
	&\quad + \f1{t^2} \int_{t/4}^{7t/4} \ju^{-2\xg} |f(r,\xo) \chi'(r/t)|^2 r^2 dr.
\end{align*}
\end{proof}

\

The following result is an analogue of Theorem 5.3 in \cite{LiT}. 

\begin{lemma} \label{inptdcExt}
Let $T$ be fixed and $\p$ solve \eqref{main} for the times $t\in [T,2T]$. 
There is a fixed positive integer $k$ such that for any multi-index $J$ with $|J| + k \le N$, we have:
\begin{align}\label{u/v decay}
\begin{split}
|\p_{J}| &\ls_{|J|} \|\p_{|J| \leq \cdot \leq |J|+k}\|_{LE^1[T, 2T]} \f{\ju^{1/2}}\jv. %
\end{split}\end{align} %
\end{lemma}

\begin{proof}We prove this by looking separately at $(t,x)$-pair values in $\crt,C^T_R$ and $\cut$. 

\begin{itemize} 

\item (The $\cut$ regions, with $1\le U \le 3T/8$)
In contrast to the ``near'' region $C^R_T$ and the ``far'' region $C^T_R$, the regions close to the cone will proceed differently: we utilise a Hardy-like inequality adapted to the cone, namely \cref{con.Hardy}. 

	Let $\chi : \R_+ \to \R_+$ be a smooth cutoff function with $\chi(s)=1, s\geq 1/2$ and $\chi(s)= 0, s\leq 1/4$.  For any smooth-enough function $w$,
\begin{align}
\begin{split}\label{u-Har}
  \|\f{w}U\|_\ltcut 
 	&\ls \Big\|\frac{\chi(\f{r}t) w}{\ju}\Big\|_{L^2[T, 2T]L^2} \\
	&\lesssim \|\partial_{r} (\chi(\f{r}t) w)\|_{L^2[T, 2T]L^2} + T\inv \| \chi(\f{r}t) w \|_{L^2_{t,x}([T,2T] \times \{ T/8 \le r  \le 15/8 T \} )}               \\
	&\ls T^{1/2} \|w\|_{LE^1[T, 2T]}
\end{split}
\end{align}where the second line follows by \eqref{con.Hardy}.

Thus
\begin{align*}
\begin{split}
\|\phi_J\|_\licut &\ls 
   \sum_{i\le 1,j\le 2} \f1{(T^3U)^{1/2}}  \|S^i \xO^j \phi_J\|_{L^2(\ti C^U_T)} +  \lr{\f{U}{T^3}}^\f12 \|\pa_r S^i \Omega^j \phi_J\|_{L^2(\ti C^U_T)} \\
 &\ls \lr{\f{U}{T^3}}^\f12 T^{1/2} \sum_{i\le1,j\le2} \| S^i \xO^j \phi_J \|_{LE^1[T,2T]} \\
 &\ls \f{U^{1/2}}T \|  \phi_{|J| \le \cdot \le |J|+k}\|_{LE^1[T,2T]}.
\end{split}
\end{align*}

\item (The $\crt$ regions, for $R$ values sufficiently small relative to $T$) This is essentially \cref{Sob.emb}: apply the Sobolev embedding estimate \eqref{in-R-Sob} to $\phi_J$
\begin{align*}
\!  \| \p_J\|_{L^\infty(\crt)} 
  &\ls \sum_{i\le 1,j\le 2}
\f1{(R^3T)^{1/2}} \|S^i \xO^j \p_J\|_{L^2(\ti C^R_T)} 
   + 
\f1{R^{1/2} T^{1/2}} \|\pa_r S^i \xO^j  \p_J\|_{L^2(\ti C^R_T)} \\
&\ls \f1{T^{1/2}} \|\phi_{|J|\le\cdot\le|J|+k}\|_{LE^1[T,2T]},
\end{align*}
and take the supremum over, say, $R< 3T/8$. The second inequality comes from commuting $S^i\xO^j$ with $Z^J$ in a way that will put it in the form \cref{vf defn}. This is where the integer $k$ arises. 

\item (The $C^T_R$ regions)
 \eqref{R-Sob} implies
\begin{align*}
\| \phi_J \|_{L^\infty(C^T_R)} &\ls
\frac{1}{R^{1/2}} \sum_{i\le1, j \leq 2}
    \|R^{-3/2}S^i \xO^j \phi_J\|_{L^2(\ti C^T_R)} +\|R^{-1/2}\pat S^i \Omega^j  \phi_J\|_{L^2(\ti C^T_R)}\\
&\ls \f1{R^{1/2}} \sum_{i\le1,j\le2} \|S^i\xO^j \phi_J\|_{LE^1[T,2T]} \\
&\ls \f1{R^{1/2}} \|\p_{|J|\le\cdot\le|J|+k}\|_{LE^1[T,2T]}.
\end{align*}
Then we take the supremum over the relevant $R$ values. In $C^T_R$, we have $v\sim r$ and $u\sim r$.%
\end{itemize}
\end{proof}

\section{Derivative estimates in $L^2$}\label{sec:der ests}

\begin{lemma}\label{dya bd} 
  Suppose that $\xs$ and $\xd$ from \cref{P def} are \emph{nonnegative} real numbers, and $\xd' \in [-1,\iy)$. Let $L,L'$ denote dyadic numbers of the form \cref{dyadic numbering}, with $L, L'=1$ when $h =0$ and, in general, $L , L' \gg_{h,g^\xo} 1$ are appropriately large relative to 1, depending on $h$ and $g^\xo$.\footnote{For example, if $h\in S^Z(\eps\jr^{-1})$ for a sufficiently small $\eps>0$, then $L=1$.}
\begin{itemize}
\item If $ L \le U,R \leq 3T/8$, then
\begin{itemize}
\item
\begin{equation}\label{crt bound}
  R\| \de w_{\le m}\|_{\ltcrt} 
  	\ls  \|w_{\le m}\|_\ltcrtt + \|Sw_{\le m}\|_\ltcrtt + R^2 \| (P w)_{\le m}\|_{L^2(\tilde C_{T}^R)}
\end{equation}
\item Let $C^U_{T,1} := \cut \cap \{ r < t\}$ and $C^U_{T,2} := \cut\cap\{r>t\}$. 
\begin{equation}\label{cut bound}
  U\| \de w_{\le m}\|_{L^2(C_{T,1}^U)} 
  \ls
 \| w_{\le m}\|_{L^2(\tilde C_{T,1}^U)}
  +  \|S w_{\le m}\|_{L^2(\tilde C_{T,1}^U)}
  + UT \| (P w)_{\le m}\|_{L^2(\tilde C_{T,1}^U)}
\end{equation}
\begin{equation}
  U\| \de w_{\le m}\|_{L^2(C_{T,2}^U)} 
  \ls
 \| w_{\le m}\|_{L^2(\tilde C_{T,2}^U)}
  + \sum_{\bar Z\in\{\xO,S\}} \|\bar Z w_{\le m}\|_{L^2(\tilde C_{T,2}^U)}
  + UT \| (P w)_{\le m}\|_{L^2(\tilde C_{T,2}^U)}
\end{equation}
\end{itemize}
\item
If $L' \le T < R$, i.e. $L' \le T \le 3R/8$, then
\begin{align}\label{ctr bound}
R\|\de w_{\le m}\|_{L^2(C^T_R)} &\ls \|w_{\le m}\|_{\lt(\ti C^T_R)} 
   + \sum_{\bar Z\in\{\xO,S\}} \|\bar Z w_{\le m}\|_{\lt(\ti C^T_R)} + R^2 \|(Pw)_{\le m}\|_{\lt(\ti C^T_R)} 
\end{align}
\end{itemize}

\end{lemma}

\begin{proof}%
We begin by proving \eqref{crt bound}. Let $w$ denote a reasonably smooth function. We shall first prove that for $1 \ll R \leq 3T/8$,
\begin{equation} \label{easy}
  R\| \de w\|_{\ltcrt} 
  	\ls  \|w\|_\ltcrtt + \|Sw\|_\ltcrtt + R^2 \| P w\|_{L^2(\tilde C_{T}^R)}
\end{equation}
	Let $\chi(t,r)$ be a radial cutoff function on $\R^{1+3}$ with $\supp\chi \subset \crtt$ and $\chi=1$ on $\crt$; a further fixing of $\chi$ will come later in the proof.  
Two observations are in order:
\begin{enumerate}
\item
If $r < t$ then for a sufficiently large constant $C'$, we have %
\begin{equation}\label{1f}
\chi \lr{ \f{u}t |\de w(t,x) |^2  }\le \chi\lr{   |\grad w |^2- w _t^2+ \f{C'}{ut} |S w |^2 }
\end{equation} (which holds without the multiplication by $\chi$ as well)
as an expansion of the terms $|Sw|^2, |\de w|^2$ reveals; the values $C' \geq 3$ work for every $(r,t)$ such that $0\le r<t$.

\item
By integration by parts,
\begin{align}\label{ini.com}
\int \chi ( |\grad w |^2 -  w _t^2) \, dxdt
	&=  \int \chi w (\pat^2-\Delta) w \,dxdt - \int \f12 (\pat^2-\Delta)\chi  w ^2\, dxdt.
\end{align} There are no boundary terms in either time or space because of the compact support of $\chi(t,r)$ in both time and space. 
\end{enumerate}
 	Integrating \cref{1f} in spacetime, we have via \cref{ini.com}
\begin{align}\label{1g}
\int\chi \f{u}t |\de w |^2 \,dxdt
	&\le \int \chi w(\pat^2-\Delta) w + O(|\Box\chi| w ^2) + \f{C'}{ut}\chi |S w |^2 \,dxdt.
\end{align}
The proof of \cref{easy} will be complete once we incorporate $Pw$ into \cref{1g}:
\begin{itemize}\label{list}
\item
Let $\Box_h$ denote the second order operator
$$\Box_h:=\paa\hab\pab.$$
For $\int (\chi w)(\Box_hw) \,dxdt$, we integrate by parts and use Cauchy-Schwarz. A term
$$\int \chi \hab \paa\p \pab \phi \,dxdt = O\lr{ \int \chi \f{ |\de w|^2}\jr \,dxdt }$$
arises, and for this term we use the hypothesis that $L\gg_h 1$ for $h \neq 0$. 

Similarly, $\int (\chi w)  (g^\xo\Delta_\xo w) \,dxdt$ is treated by integration by parts and Cauchy-Schwarz. We use the smallness of $\jr^{-2-\xd'}$ (which is $O(\jr\inv)$ since $\xd'\in[-1,\iy)$) for sufficiently large $R$.
\item
We use the bound $V\ls \jr^{-2}$.

\item
For $\int \chi w B \pa w$ we use Cauchy-Schwarz. For $\int \chi w \pa(Aw)$ we integrate by parts and use Cauchy-Schwarz; it is also possible to bound this using information on $\pa A$ if one does not integrate by parts, but we integrate by parts in order to use fewer assumptions. The bounds we obtain are sufficient to prove the claim \cref{easy} even when $\sigma = 0$, and we only assume $A,B\in S^Z(\jr^{-1})$ in this part.
\end{itemize}

Assuming $\Box\chi\ls \jr^{-2}$, separating
$|\chi wPw| \ls \chi [ (R\inv w)^2 + (RPw)^2 ]$ in the right-hand side of 
\cref{1g}, and using the reasoning in the bullet points (along with the triangle inequality) to deal with $\int (\chi w) ( (\Box-P) w) \,dxdt$,
this proves the claim \cref{easy} for $C^R_T$.

\

The same proof shows the analogue of \cref{easy} for the $C^U_T \cap \{ r < t\}$ region, 
\begin{equation} \label{easy1}
  U\| \de w\|_{L^2(C^U_T\cap \{ r < t\} )} 
  	\ls  \|w\|_{L^2(\ti C^U_T \cap \{ r < t\})}  + \|Sw\|_{L^2(\ti C^U_T \cap \{ r < t\})} + UT \| P w\|_{L^2(\ti C^U_T \cap \{ r < t\})}
\end{equation}
if we choose a $\chi$ adapted to $\cut \cap \{ r < t\}$ (rather than $\crt$) that satisfies $\bo\chi \ls \f1{\la t+r\ra\nm} $ (rather than $\bo\chi \ls 1/\jr^2$).\footnote{(Note that if $T$ is sufficiently large, then we may even take $L = 1$ for $\cut$ and $L'=1$ for $C^T_R$.)}

\

Similar arguments show the result for vector fields, \cref{crt bound} and \cref{cut bound}. The only new thing one has to deal with is $\int \chi w_{\le m} [P,Z^{\le m}] w \,dxdt$ and similar arguments involving integration by parts and Cauchy-Schwarz establish the claims \cref{crt bound,cut bound}. 

\

Next, we prove 
\begin{equation} \label{eas1}
R\|\de w\|_{L^2(C^T_R)} \ls \|w\|_{\lt(\ti C^T_R)} 
   + \sum_{\bar Z\in\{\xO,S\}} \|\bar Z w\|_{\lt(\ti C^T_R)} + R^2 \|Pw\|_{\lt(\ti C^T_R)}.
\end{equation}
	The proof for the region $\{r > t\}$ is essentially a switching of the $r$ and $t$ variables in 
what has been done 
for the $C^R_T$ and $C^U_T\cap \{ r < t\}$ regions. 
	For any point $(t,x)$ such that $|x| > t$, 
\begin{equation}\label{honest}
|\de w(t,x)|^2 \le \f{r}{r-t}(w_t^2 -  w_r^2) + \f{C'}{(r-t)^2}(Sw)^2 + C \f{(\xO w)^2}{r^2}
\end{equation} 
for some sufficiently large constants $C,C'>0$. For the angular derivatives, this follows because $\pao = \sum_j c_j \xO_j$ for some coefficients $c_j$ such that $|c_j| \ls 1/r$. %
We shall only use the weaker estimate 
\begin{equation}\label{wea}
|\de w(t,x)|^2 \le \f{r}{r-t}(w_t^2 -  |\grad w|^2) + \f{C'}{(r-t)^2}(Sw)^2 + C \f{(\xO w)^2}{r(r-t)}. 
\end{equation}
We use this because it makes \eqref{dis} conceptually cleaner; and because using \eqref{honest} would lead to no gain in the final derivative estimates for $C^T_R$, due to the presence of the $(r-t)^{-2}$ coefficient of $(Sw)^2$.

\begin{itemize}
\item (Bound in $C^T_R$) Let $\chi(t,r)$ be a radial cutoff function adapted to $C^T_R$. By \cref{wea},
\begin{align}\label{dis}
\begin{split}
\int\chi|\de w |^2 \,dxdt
	&\leq \int \f{r}{r-t} \chi(w_t^2 - |\grad w|^2) + C\f{r}{r-t}\chi |\f{(\xO w)^2}{r(r-t)}| + \f{C'} {(t-r)^2} \chi |S w |^2 \,dxdt. %
\end{split}
\end{align} %
The analysis henceforth is similar to the three bullet points above.
Assuming $\Box\chi \ls \jr^{-2}$, 
	we end up with 
$$\|\de w\|_\ltctr \ls R\inv \lr{ \|w\|_\ltctrt + \sum_{\bar Z\in\{\xO,S\}} \|\bar Z w\|_\ltctrt } + R\|Pw\|_\ltctrt,$$
i.e., \cref{eas1}. 

\item (Bound in $C^U_T\cap \{ r > t\}$)
We adapt $\chi$ to $C^U_T \cap \{ r > t\}$ with $\Box\chi \ls (\la t+r\ra\nm)\inv$. %
Then by Cauchy-Schwarz,
\[ \|\de w\|_{\lt(\cut \cap \{ r > t\})} 
\ls  U\inv \lr{ \|w\|_{\lt(\cutt\cap \{ r > t\})} + \sum_{\bar Z\in\{\xO,S\}} \|\bar Z w\|_{\lt(\cutt\cap \{ r > t\})} } 
     + T \|Pw\|_{\lt(\cutt\cap \{ r > t\})}. \] %
\end{itemize}

The full results for vector fields $w_{\le m}$ again follow simply by similar integration by parts and Cauchy-Schwarz arguments.
\end{proof}

We will need to bound the second derivative of vector fields in $L^2$ when proving $L^\iy$ estimates for vector fields of a function. Hence we present \cref{dya bd 2} immediately.

\begin{corollary}\label{dya bd 2}
Assume the hypotheses of \cref{dya bd}. Then
\begin{itemize}
\item
\begin{equation}
  R\| \de^2 w_{\le m}\|_{\ltcrt} 
  	\ls  \|\de w_{\le m+n}\|_\ltcrtt + R^2 \| \de (P w)_{\le m}\|_{L^2(\tilde C_{T}^R)}
\end{equation}
\item 
\begin{equation}
  U\| \de^2 w_{\le m}\|_{L^2(C_{T}^U)} 
  \ls
 \| \de w_{\le m+n}\|_{L^2(\tilde C_{T}^U)}
    + UT \| \de (P w)_{\le m}\|_{L^2(\tilde C_{T}^U)}
\end{equation}
\item
\begin{equation}
R\|\de^2 w_{\le m}\|_{L^2(C^T_R)} 
\ls 
 \|\de w_{\le m+n}\|_{\lt(\ti C^T_R)} + R^2 \|\de (Pw)_{\le m}\|_{\lt(\ti C^T_R)} 
   \end{equation}
\end{itemize}
\end{corollary}

\begin{proof}
Fixing any $\alpha\in\{0,1,2,3\}$ and denoting $\paa$ by $\pa$,
 we substitute $\pa w_{\le m}$ for the function $w$ in the proof of \cref{dya bd}. 
 
A new type of term arises, which is 
$$\int \chi \pa w_{\le m} P\pa w_{\le m} = 
\int \chi \pa w_{\le m} ( \pa f_{\le m} + \pa[P,Z^{\le m}]w + [P,\pa] w_{\le m}).$$
We can handle the first term on the right-hand side by Cauchy-Schwarz. 

For the $\Box_h,g^\xo\Delta_\xo$ and $V$ contributions to $P$, similar arguments as before using Cauchy-Schwarz and integration by parts work. For the contributions of the $\paa A^\x$ and $B^\x\paa$ components to $P$ in both $\pa [P,Z^{\le m}]w$ and $[P,\pa] w_{\le m}$, we also use integration by parts and Cauchy-Schwarz, and the fact that $\pa A \in S^Z(\jr^{-2})$, 	i.e. this bound holds for all $(r,t)$ (and hence all three dyadic regions),
	which follows from the assumptions $\pa A\in S^Z_{int}(\jr^{-2}) \cap S^Z_{cone}(\jr^{-2})$ and $\pat A\in S^Z(\jv\ju\inv \jr\inv \jr^{-1-\sigma})$ because of \cref{D decomp}. 
	More concretely, we have the schematic equalities
$$\int \chi \pa w_{\le m} \pa[B^\x\paa, Z^{\le m}] w = \int (\chi' \pa w_{\le m} + \chi \pa^2w_{\le m}) \ti B \pa w_{\le m}$$
$$\int \chi \pa w_{\le m} \pa[\paa A^\x, Z^{\le m}] w = \int (\chi' \pa w_{\le m} + \chi \pa^2 w_{\le m}) \pa \ti A \cdot w_{\le m} $$
where tildes denote vector fields. We apply the aforementioned assumptions $\pa A\in S^Z(\jr^{-2})$ and $B\in S^Z(\jr^{-1})$. 

\end{proof}
 
\begin{corollary}[$L^\iy$ estimates for derivatives] \label{derdecbtr} 
Assume the hypotheses of \cref{dya bd 2}. Hence $\xs$ and $\xd$ from \cref{P def} are \emph{nonnegative} real numbers. 
\begin{enumerate}
\item If $1 \ll U \leq 3T/8$, we have
\begin{align}\label{first infty}
\begin{split}
\|\pa w_{\le m}\|_\licut 
	&\ls \f1{\sqrt{UT^3}} \Big( U\inv\|w_{\le m+n}\|_\ltcutt + T (\|(Pw)_{\le m+n}\|_\ltcutt + \|U \pa(Pw)_{\le m}\|_\ltcutt ) \Big).
\end{split} 
\end{align}

\item Let $1 \ll R \leq 3T/8$. Then we have: 

\[ \|\pa w_{\le m}\|_\licrt \ls \f1{\sqrt{TR^3}} \lr{ 
   R\inv \|w_{\le m+n} \|_\ltcrtt + R ( \|(Pw)_{\le m+n}\|_\ltcrtt + \|R \pa(Pw)_{\le m}\|_\ltcrtt ) 
    }	.
\]
\item Let $1 \ll T \leq 3R/8$. Then we have: 

\[ \|\pa w_{\le m}\|_\lictr \ls \f1{\sqrt{TR^3}} \lr{ 
   R\inv \|w_{\le m+n} \|_\ltctrt + R ( \|(Pw)_{\le m+n}\|_\ltctrt + \|R \pa(Pw)_{\le m}\|_\ltctrt ) 
    }	.
\]\end{enumerate}
\end{corollary}

\begin{proof} Let $v = w_{\leq m}$. The main idea in this proof is to 
\begin{itemize}
\item
first use the initial $L^\iy$ estimates proved in \cref{sec:initialests} on \emph{derivatives} $\pa v$, and to commute this $\pa$ with the vector fields $S^i\xO^j$ in both terms of the majorizer in the estimates \cref{U-Sob,R-Sob,in-R-Sob}. This results in
\begin{align*}
\|\pa v\|_\iy  &\ls \sum_{i\le 1,j\le 2} (W^3W')^{-1/2} \|S^i\xO^j \pa v\|_2 + (\ti W)((W^3W')^{-1/2}) \| \de S^i\xO^j \pa v\|_2  \\
  &\ls (W^3W')^{-1/2} \|\pa v_{\le n}\|_2 + (\ti W)((W^3W')^{-1/2}) \|\pa^2 v_{\le n}\|_2 \\
  &= (W^3W')^{-1/2} \Big( \|\pa v_{\le n}\|_2 + \ti W \|\pa^2 v_{\le n}\|_2  \Big)
\end{align*} 
for dyadic weights $W,W'$ and $\ti W\in \{W,W'\}$, where the choices of $W,W'$ and $\ti W$ all depend on the region in question.
\item
And secondly to use the derivative estimates just proved in \cref{dya bd} and \cref{dya bd 2}, in order to control $\|\de v_{\le n}\|_2$ and $\ti W \|\de^2 v_{\leq n}\|_2$ respectively. 
\end{itemize}

In $\cut$, one has $W = T$ and $W' = \ti W = U$.  Let $k\ge0$ be any integer. Then
\begin{align*}
\| \pa v_{\le k} \|_2 + \ti W \|\pa^2 v_{\le k}\|_2
  &= \| \pa v_{\le k} \|_2 + U \|\pa^2 v_{\le k}\|_2 \\
  &\ls \|\pa v_{\le k+1}\|_2  + UT \| \pa (Pv)_{\le k}\|_2 \\
  &\ls U^{-1} \|v_{\le k+2}\|_2 + T \|(Pv)_{\le k+1}\|_2 + UT \| \pa (Pv)_{\le k}\|_2
\end{align*}
This proves \cref{first infty}. For the other two regions, the proof is similar.
\end{proof}

\section{Setup for pointwise decay iteration}\label{sec:setup}

\begin{definition}
Let 
\begin{equation*}
\rei := \{ \text{dyadic numbers }R : R \ge1, R < \f{t-r}8 \}
\end{equation*}
denote the collection of dyadic numbers we call Region 1, and let \begin{equation*}
\reii := \{ \text{dyadic numbers }R : R \ge1,  \f{t-r}8 \le R < t+r \} %
\end{equation*}
denote the collection we call Region 2.
\end{definition}

\begin{definition}
Let $\R_+ :=[0,\iy)$. 
\begin{itemize}
\item
Let $D_{tr}$ denote   
	\[
	D_{tr} := \{ (\rho,s)  \in \R_+^2: -(t+r) \leq s-\rho\leq t-r, \ |t-r| \leq s+\rho \leq t+r\}.
	\] 
When we work with $D_{tr}$ we shall use $(\rho,s)$ as variables, and $D_{tr}^{ R}$ is short for $D_{tr}^{\rho\sim R}$.
\item
For $R>1$, let 
\begin{equation*}
\co:= D_{tr} \cap \{ (\rho,s) : R < \rho<2R\}
\end{equation*}
and let 
\[
D^{R=1}_{tr} := \dtr\cap\{(\rho,s) : \rho<2 \}. \]
\end{itemize}
\end{definition}

\begin{lemma}[Maximal vertical length within $D_{tr} \cap \{ (\rho,s)\in \R_+^2 : \rho \le s\}$] \label{height bound}
Uniformly in the set of $r,t$ values lying in $\{ (r,t) : 0 \le r \le t\}$, we have that for any point $(\rho',s')\in\dtr \subset \R^+_{\rho'} \times \R^+_{s'}$,
\begin{enumerate}
\item
If $r \le t/3$, then
\begin{equation*}
|\dtr\cap \{(\rho',s') : \rho = \rho'\} | \le \min\{2\rho,2r\}
\end{equation*}
\item
If $t\ge r\ge t/3$, then
\begin{equation*}
|\{ s'\ge \rho'\ge0\} \cap \dtr\cap \{(\rho',s') : \rho = \rho'\} | \le t-r
\end{equation*}
\end{enumerate}
where $| \cdot|$ denotes the length.
\end{lemma}
\begin{proof} We split the proof into two cases. 
\begin{enumerate}
\item
Let $r \le t/3$; then for each $\rho$, the maximal vertical length within $\dtr$ is $2r$ and occurs when $r \le \rho \le \f{t-r}2$; by symmetry, this length, $2r$, is maximal. When $0 \le \rho \le r$, the maximal vertical length of $\dtr$ is $2\rho$, which implies that this value of this length is sharp if and only if $0 \le \rho \le r.$

\item
Let $r \ge t/3$; then for each $\rho$, the maximal vertical length within $\dtr \cap \{ s \ge \rho \}$ is $t-r$ and occurs when $\f{t-r}2 \le \rho \le r$ and by symmetry once more, this length, $t-r$, is maximal. Furthermore, in a manner precisely analogous to the $r\le t/3$ case, we once more have that when $0 \le \rho \le \f{t-r}2$, the bound $2\rho$ is sharp if and only if $\rho$ lies in this small region.
\end{enumerate}
\end{proof}

\begin{definition} Given $\xl\in\R$,
\[ \xk(\xl,t-r) := \begin{cases}
1 & \xl>1\\
\log\nm & \xl=1\\
\nm^{1-\xl} & \xl<1
\end{cases}. \]
In this paper, this function arises either as 
\[\sum_{1\le R\ls\nm} \f1{R^{\xl-1}} \ \ \text{ or  }\ \ \int_0^{t-r} \f1{\jxu^{\xl}}d\xu.\]
\end{definition}

\begin{lemma}\label{conversion}
Let $m\ge0$ be an integer and suppose that $\psi : [0,\iy)\times\R^3\to\R$ solves $$\Box\psi (t,x)= g(t,x), \qquad (t>0, \ x\in\R^3) $$ with vanishing initial data, with 
\[ g \ls \f{\log^m\nm}{ \jr^\x \jt^\xb \nm^\eta }, \]
where the values of $\x,\xb,\eta$ will be specified below. 

\begin{itemize}
\item (The case $r\le t$) 
Assume that $\xb\ge0$ and $\eta\in\R$. Assume also that $|x| \leq t$. Assume that $\supp g$ is contained within $\{ r\le t\}$. 

If $1<\x<3$, then 
\begin{align}\label{eq:conversion lemma bound}
\begin{split}
\f{\jr\psi(t,x)}{\log^m \nm}  \ls\min&\lr{ \f{\xk(\x-1,t-r)}{\nm^{\xb+\eta-1}} , \f1{\nm^{\xb+\eta+\x -3}}  }  +  \f{\xk(\eta,t-r)}{\nm^{\x+\xb-2}}
\end{split}
\end{align}

If $\x>3$ (we will not be needing the cases $\x=3$ or $\x\le1$), 
\begin{align}\label{eq:conversion lemma bound 2}
\begin{split}
\f{\jr\psi(t,x)}{\log^m \nm}  \ls \f{\xk(\x-1,t-r)}{\nm^{\xb+\eta-1}} + \f{\xk(\eta,t-r)}{\nm^{\x+\xb-2}} 
\end{split}
\end{align}

\item (The case $r>t+1$)
Let $\x > 1, \eta\in\R$. Suppose that $r>t + 1$, %
and 
$$g\ls \f{1}{\jr^\x\nm^\eta}.$$
Then 
\begin{equation}\label{r>t bd}
\jr\psi \ls 
\f1{\nm^{\x-2}}
\begin{cases}
1/\nm^{\eta-1} & \eta>1\\
1/\la r+t\ra^{\eta-1} &\eta<1
\end{cases}
\end{equation}
For the case $r>t$: in this paper we will only need $\x>1, \eta \ne 1$; a full explanation is given in the proof of \cref{exteprop}. 
\end{itemize}
\end{lemma}

\begin{proof}
\begin{enumerate}
\item (The case $r\le t$)
We write 
\[ \int_\dtr \rho \sup_{S^2} |\Box\psi|dsd\rho = \sum_\rei \int_\co \rho \sup_{S^2} |\Box\psi|dsd\rho + \sum_\reii \int_\co \rho \sup_{S^2}|\Box\psi|dsd\rho \]
and bound $\Box\psi$ pointwise by the bound in the hypotheses. Throughout $\dtr$, we have \[ \f1s \ls \f1{t-r} \] and we will use this repeatedly below.

We begin with the first bound in \eqref{eq:conversion lemma bound}, namely, 
\[
\sum_\rei \int_\co \rho \sup_{S^2} |\Box\psi|dsd\rho 
	\ls \f{\xk(\x-1,t-r)}{\nm^{\xb+\eta-1}}. 
\]
In the region $\rei := \{ 1 \leq R < \f{t-r}8 \}$ defined at the beginning of this section, we have
	\[ s-\rho \sim t-r. \]
Therefore, for $R\in \rei$ and any $\xb\ge0, \eta\in\R$,
  \begin{align*}
   \int_\co \rho \sup_{S^2} |\bo\psi|dsd\rho 
   	&\ls \int_\co |\f{\log^m \jxu}{\jrho^{\alpha-1} \js^\beta \jxu^\eta} | dsd\rho\\
   	&\ls \f{\log^m \nm}{R^{\alpha-1} \nm^\eta} \int \f{dsd\rho}{\js^\beta} \\
	&\sim \f{\log^m \nm}{R^{\alpha-2} \nm^{\eta}} \int \f{ds}{\js^\xb}   \\
	&\ls \f{\log^m \nm}{R^{\alpha-2} \nm^{\eta}} \f1{\nm^{\xb-1}}  \\
	&= \f{\log^m \nm}{R^{\alpha-2} \nm^{\beta+\eta-1}} 
  \end{align*}
where $\xu :=s-\rho$. 
Thus 
   \begin{equation*}
     \sum_{1 \leq R < \f{t-r}8 } \f{\log^m \nm}{R^{\alpha-2} \nm^{\beta+\eta-1}} 
     	= \f{\log^m \nm}{\nm^{\beta+\eta-1} } \xk(\x-1,t-r).
   \end{equation*}

Next, we prove that when $\x<3$, and $\xb\ge0, \eta \in \R$, we have
\[\sum_\rei\int_\co\rho|\bo\psi|\,dsd\rho 
	\ls  \f{\log^m\nm}{\nm^{\xb+\eta+\x-3}}.
\]
This is shown as follows: since $\xb\ge0$, we have $\js^{-\xb}\ls\nm^{-\xb}$, and
\begin{align*}
\log^{-m}\nm \sum_\rei \int_\co \rho |\bo\psi|dsd\rho 
	&\ls \nm^{-\xb-\eta} \sum_\rei R^{1-\x} \int \int  dsd\rho \\
	&\ls \nm^{-\xb-\eta} \sum_\rei R^{3-\x} \\
	&\ls \nm^{-\xb-\eta+3-\x}
\end{align*}
where the last line follows by the hypothesis $\x<3$. 

Finally, we show that when $\x>1$ and $\xb\ge0$, then
\[\int_{\bigcup_{R\in\reii}\co} \rho|\bo\psi|dsd\rho \ls \log^m\nm \f{\xk(\eta,t-r) }{ \nm^{\x+\xb-2} } \]
which will complete the proof. 
For $R\in\reii$, we employ the fact that when $\xb\ge0$ we have
\[
\jrho^{-\xb}\ls \nm^{-\xb} \]
to find that
\begin{align*}
\log^{-m}\nm \int_\co \rho |\bo\psi|dsd\rho 
	&\ls \nm^{1-\x} \int_\co \js^{-\xb}\jxu^{-\eta}dsd\rho \\
	&\ls \nm^{1-\x }\nm^{-\xb} \int ds \int_0^{t-r} \jxu^{-\eta} d\xu \\
	&\ls \f1{\nm^{\xb+\alpha-2}} \xk(\eta,t-r)
\end{align*}
with the last line following by \cref{height bound}. 
\item (The case $r>t$) We now prove \eqref{r>t bd}.  Assume that $\x>1$. A straightforward integration shows that 
\begin{align*}
\int_\dtr \rho \f{dsd\rho}{\jrho^\x \on^\eta} 
	&\ls
	\f1{\nm^{\x-2}}
	\begin{cases}
	\f1{\nm^{\eta-1}} & \eta>1\\
	\ln \f{\la r+t\ra}{\nm} & \eta=1 \\
	\f1{\la r+t\ra^{\eta-1}}& \eta<1
	\end{cases}
\end{align*}which shows \eqref{r>t bd}. 
\end{enumerate}
\end{proof}

Next, in \cref{rem:int} we apply \cref{conversion} in a setting that is relevant to our current problem:
\begin{remark}[$P$'s coefficients in the wave zone determine final decay rate]\label{rem:int}
In this remark we show how the pointwise decay rate for the solution (and its vector fields) improves indefinitely in the subset of $\dtr$ closest to $\{\rho=0\}$. We also show that it is the pointwise decay rate of the coefficients in the wave zone (i.e., close to $\{ \rho = s\}$) that determines the final decay rate of the solution. 

Let $\pmn$ satisfy the bounds
\begin{equation}
\label{starting}\pmn(t,x) \ls \f1{\jt\ju^{\x}}, \quad \x\in\R.
\end{equation}
Then the contribution of the integral over $D_{tr}^{\ll u}:=\dtr \cap \{ (\rho,s) : \rho \ll u\}$ is bounded by $\ju^{-\x-\xs}$; thus this integral always gains $\ju^{-\xs}$ (no matter what the value of $\alpha$ is), in contrast to the integration that takes place over the wave zone. In the rest of this remark, we will show this.

Suppose that 
\begin{equation}\label{the form}
\Box\p_{1} = V_{1}(t,x)\pmn
\end{equation}where $V_{1}\in S^{Z}(\jr^{-2-\xs}).$
Note that by the procedure outlined in \cref{pr}, essentially the entire problem \cref{eq:problem} can be written in this form \cref{the form}.\footnote{The only exception is an extra consideration near the light cone $\{r=t\}$, which is explained in \cref{sec:cone bounds}.} 
Then by \cref{starting}, if given a function $G$ we let
$$H_G := \sum_{k=0}^2 \|\Omega^k (G)_{\leq m+n} (t, r\omega)\|_{L^2(\bbS^2)}$$
then
\begin{align*}
\int_{D_{tr}^{\ll u} \cap \{ R <\rho<2R \}} \rho H_{\Box\p_1} \, dA 
 &\ls \int_{D_{tr}^{\ll u} \cap \{ R <\rho<2R \}} \js\inv\la s-\rho\ra^{-\alpha} \rho H_{V_1} \, dA &\textrm{by \cref{starting}}\\
 &\ls \ju\inv \int_{D_{tr}^{\ll u} \cap \{ R <\rho<2R \}} \la s-\rho\ra^{-\alpha} \rho H_{V_1} \, dA \\
 &\sim \ju^{-1-\x} \int_{D_{tr}^{\ll u} \cap \{ R <\rho<2R \}} \rho H_{V_1} \, dA \\
 &\ls \ju^{-1-\x} R^{1 - \xs}.
\end{align*}For the lowest $R$ value we integrate over $D_{tr}^{\ll u} \cap \{ \rho \ls 1 \}$. 
Summing in $\{ R : R\ll u\}$, we obtain for all sufficiently small $\xs$, namely $\xs<1$ (or, if $\xs$ is large, we can truncate to values less than 1), the following bound
\begin{align*}
\ju^{-\xs-\x}.
\end{align*}
In conclusion, we have indefinite improvement in the ``deep interior'' $\{\rho \ll u\}$ of $\dtr$ of the pointwise decay rate of $\p$, no matter the starting decay rate (i.e. no matter the value of $\x$). 

Let us now bound the integral in the wave zone. If $\alpha>1$, then by \cref{conversion}
\begin{align*}
\int_{\dtr \setminus D_{tr}^{\ll u}} \rho H_{\Box\p_1} \, dA 
  &\ls \ju^{-1-\xs}
\end{align*}
Thus, once $\x>1$, it is $V_1$'s decay rate near the cone that dictates the decay rate for $\phi_1$ of 
$$\jr\p_1 \ls \ju^{-1-\xs}.$$ 
Using a later tool (\cref{r-t}), this bound can be improved to 
$$\p_1\ls\jv\inv\ju^{-1-\xs}.$$

\end{remark}

\begin{definition}[Cutoff functions]\label{cutdef}
Let $\chi_\text{exte}(t,x)$ denote a smooth radial cutoff function adapted to $\{r \ge t,  r-t \sim r \}$.

Let $\cinte(t,x)$ denote a smooth radial cutoff function adapted to $\{ r \leq t, t-r \sim t\}$. 

Let $\cone(t,x)$ be a smooth radial cutoff function equalling $1 - (\cinte +\cexte)$.
	We also assume $\supp\cone\subset\{r/2 \leq t \leq 3t/2\}$. 

Thus: in $C_T$, for instance, $\cinte$ and $\cone$ sum to 1, while in $([T,2T] \times \rt)\setminus C_T$, $\cexte$ and $\cone$ sum to 1. 
\end{definition}

In the following sections, we shall finish the proof of \cref{thm:main}; by the product rule, and also \eqref{D decomp}, it will suffice to prove pointwise decay for
\begin{align}\label{pr}
\begin{split}
\Box\pm &= (\ti V + \pa \ti B)\pm + \pa( [\ti A + \ti B]\pm) + \pa(\ti h\pa\pm) + \ti {g^\xo} \pa^2 \pm, \\
(\pm(0), \bar N\pm(0)) &= ( (\phi_0)_{\le m}, (\phi_1)_{\le m-1} ) 
\end{split}
\end{align} 
where $\,\ti \hfill\,$ denotes vector fields. 

Before commencing the pointwise decay iteration in the next section, we note that:
\begin{itemize}
\item
By \cref{u/v decay}, the desired decay rate in \cref{thm:main} already holds in the region $\{ |u|\le1\}$. Henceforth in this article, we shall assume that $|u|>1$, i.e., $|t-r|>1$. Thus we work away from the light cone $\{r=t\}$. 
\item
Due to the domain of dependence properties of the wave equation, we shall first complete the iteration in $\{r>t+1\}$, which is the content of \cref{sec:exte}. For the iteration in $\{r < t-1\}$, the decay rates obtained from the fundamental solution are insufficient in the region $\{r < t/2\}$. To remedy this, we prove \cref{r-t}. With the new decay rates obtained from \cref{r-t}, we are then able to obtain new decay rates for the solution and its vector fields. At every step of the iteration, \cref{conversion} is used to turn the decay gained at previous steps into new decay rates.
\end{itemize}

\section{The upper bound in $\{ r>t+1\}$}\label{sec:exte}

Before embarking on the pointwise decay iteration for the equation in \cref{eq:problem}, we first explain in \cref{rem:id} how we deal with the initial data in \cref{eq:problem}.

\begin{remark}[The initial data] \label{rem:id}
Let $w := S(t,0)(\p_0,\p_1)$ denote the solution to the free wave equation with initial data $(\p_0,\p_1)$ at time 0. Thus 
$$w_J(t,x) = \f1{|\pa B(x,t)|} \int_{\pa B(x,t)} (\p_0)_J(y) + \nabla_y (\p_0)_J(y) \cdot (y-x) +t (\p_1)_J(y) \, dS(y).$$
Let $\x -1\in \{1+\min(\xs,\xd,\xd'), 1+\min(\xs+1,\xd,\xd')\}$ with the first (resp. second) one as the value of $\x-1$ assuming hypotheses from part 1 (resp. part 2) of \cref{thm:main}. For any multiindex $J$, we now show that
\[ w_J \ls \f1{\jv\ju^{\x-1}}\] 
by the Kirchhoff formula and the weighted $\lt$ decay assumption on the initial data. We use Cauchy-Schwarz and Sobolev embedding to control the free wave pointwise by the weighted $L^2$ bound assumed on the initial data.
When $r\gg t$ and $y\in \pa B(x,t)$, 
$$|(\p_0)_J(y)| + |\nabla (\p_0)_J(y)\cdot (y-x)| + |t(\p_1)_J(y)|
	\ls \jr^{-\x}$$
so that $w_J \ls \jr^{-\x} \ls \jv\inv\ju^{-(\x-1)}$.	
	Similarly, when $r \ll t$ and $y\in \pa B(x,t)$, 
$$|(\p_0)_J(y)| + |\nabla (\p_0)_J(y)\cdot (y-x)| + |t(\p_1)_J(y)|
	\ls \jt^{-\x}$$
so that $w_J \ls \jt^{-\x} \ls \jv\inv\ju^{-(\x-1)}.$ When $r\sim t$, we have $w_J \ls \jv\inv$. 
\end{remark}

\

Recalling \cref{pr}, in this section we prove that the solution to
\begin{equation}\label{wm eqn}
\Box w_{(m)} = O(\jr^{-2-\min(\xd', \xd, \xs + 1)} )\pmn
\end{equation} 
obeys the full rate bound 
$$\wm \ls \f1{\jr\nm^{1+\min(1+\xs,\xd, \xd')}}$$ 
in $\{ r>t+1\}$
assuming that $ B\in S^Z(\jr^{-2-\xs}), \pat B\in S^Z(\jr^{-3-\xs})$. 
	We used the results from \cref{sec:der ests} in transitioning from \cref{pr} to \cref{wm eqn}. 

If $ B\in S^Z(\jr^{-1-\xs}), \pat B\in S^Z(\jr^{-2-\xs})$, then we instead have $$\Box w_{(m)} = O(\jr^{-2-\min(\xd', \xd, \xs)} )\pmn$$ and the final bound $$\wm \ls \f1{\jr\nm^{1+\min(\xs,\xd, \xd')}};$$ the argument shown below proving \cref{exteprop} covers this case equally well. For the sake of simplicity and concreteness, we pick and fix the assumption \eqref{wm eqn}. 

\eqref{wm eqn} includes all the terms in \eqref{pr} except for the parts of the right-hand side of \eqref{pr} that are supported near the cone; we prove estimates for those parts in \cref{sec:cone bounds}. 

\begin{proposition}\label{exteprop}
Assume that $r>t+1$. Assuming the hypotheses of part 2 of \cref{thm:main}, 
$$\wm \ls \f1{\jr\ju^{1+\min(1+\xs,\xd,\xd')}}.$$
Assuming the hypotheses of part 1 of \cref{thm:main}, 
$$\wm \ls \f1{\jr\ju^{1+\min(\xs,\xd,\xd')}}.$$
\end{proposition}

\begin{proof} We only prove the case assuming the hypotheses of part 2 (thus $\Box \wm = O(\jr^{-2-\min(\xd', \xd,\xs+1)}) \pmn$) since the other case is similar. 
  
Given $\Box \wm = \bar G \pmn$ with $\bar G = O(1/\jr^{1+\beta})$ (here, $\beta = 1+\min(\xs+1,\xd,\xd')$), the first step is to use \eqref{u/v decay} and \cref{conversion}, which yields
\begin{equation}\label{firs}
\wm \ls \f{1}{\jv^{-1/2}\ju^\xb}.
\end{equation}
Then, a second application of \cref{conversion} yields
$$\wm \ls \begin{cases}
\f1{\jv\ju^{2\xb-5/2}} & \xb-\f12>1 \text{ i.e. } \min(\xd',\xd,\xs+1) > \f12 \\
\f1{\jv^{\xb-1/2}\ju^{\xb-1}} & \xb-\f12<1 \text{ i.e. } \min(\xd',\xd,\xs+1) < \f12
\end{cases}.
$$
Note that the sum of exponents in the denominator, call it $i_n$ if we are at step $n$, has increased by $\min(\xs+1,\xd,\xd')$. 

The case $\eta = 1$ in \cref{conversion}, whenever $r>t$, arises if $n\min(\xd,\xd') = 1$ for some integer $n\ge1$, but in this case we incur an arbitrarily small polynomial loss in $\nm$; so we avoid having to apply \cref{conversion} for the case when $\eta$ takes the value 1. For $r>t+1$, given a certain fixed value of $i_n$, we always have $i_{n+1} - i_n = \min(\xs+1,\xd,\xd')$ or $i_{n+1} - i_n = \min(\xs+1,\xd,\xd')-$, with the latter occurring if there is a borderline case which leads to an arbitrarily small loss. 

Let $a:= \min(1+\xs,\xd,\xd')$.  The general pattern after the first iterate \eqref{firs} %
is as follows. Suppose that $a = \min(\xd,\xd') <1/2$ (the case $a=1/2$ is similar because we just incur an arbitrarily small polynomial loss). 		
For some integer $N\ge1$, one has either
\begin{equation}\label{weak bds}
\text{(A)  }w_{(m)}\ls \jr^{-\f12} \jr^{-N\ti a} \ju^{-N\ti a} \qquad \text{or} \qquad\text{(B)  }w_{(m)}\ls \ju^{-\f12} \jr^{-N\ti a} \ju^{-(N+1)\ti a}
\end{equation}
for some $\ti a \in (0,a]$ which we may (and do) choose to be arbitrarily close to $a$ in the event of a borderline case, whereas $\ti a=a$ if and only if we are in a non-borderline case. %

The pattern will cycle between these two, starting at (A) for an integer value $N$, going to (B) for that same integer $N$, and then going to (A) for the integer value $N+1$, then to (B) for the integer $N+1$, and so on. 

There are correspondingly two kinds of integrals, as follows: writing $N$ now for the previous value of $N+1$, so that we work with $a$ instead of $\ti a$, by \cref{conversion} the function $rw_{(m)}$ is bounded by one of
\begin{align}\label{weak ints}
\begin{split}
\min \lr{ \f1{\ju^{\f12+(N+1)a}} , \f1{\ju^{1+a}} } \cdot
  \begin{cases}
     \f1{\ju^{Na-1}} & \text{if } Na>1\\
     \f1{\jr^{Na-1}} & \text{if } Na<1
     \end{cases}, 
\ \ 
\f1{\ju^{(N+1)a}}  \begin{cases}
     \f1{\ju^{(N+1)a-\f12}} & \text{if } \f12 + (N+1)a>1\\
     \f1{\jr^{(N+1)a-\f12}} & \text{if } \f12 + (N+1)a<1
     \end{cases} 
\end{split}
\end{align}
This iteration continues until $Na>1$. Then respectively
$$\wm \ls \f1{\jr\nm^{-\f12 + (2N+1)a}}, \qquad 
	\wm \ls \f1{\jr\nm^{-\f12+2(N+1)a}}.$$
For the %
minimal integer $N$ satisfying $Na>1$, by \cref{conversion} we have 
$$\jr\wm\ls 1/\nm^{1+a - \f32 + (2N+1)a} \le 1/\nm^{1+a}, \qquad \jr \wm\ls 1/\nm^{1+a -\f32 + 2(N+1)a} \le 1/\nm^{1+a}.$$

Suppose that $a > 1/2$. After one iteration, by \cref{conversion}, $w_{(m)}\ls \jv^{3/2}/\ju^{1+\min(1+\xs,\xd,\xd')}$. After the second iteration, $$w_{(m)}\ls \f1{\jr\ju^{-\f12+2a}}.$$ In the third iteration, one obtains $1/\ju^{1+\min(1+\xs,\xd,\xd')}$ from the $\rho$ integration alone, and for $a$ which is big enough ($a > 3/4$,
more precisely), the iteration halts here. For $1/2 < a \le 3/4$, continuing as many times as necessary, one eventually gets $$w_{(m)}\ls \f1{\jr\ju^{1+\min(1+\xs,\xd,\xd')}}.$$
\end{proof}

\section{Converting $r$ decay to $t+r$ decay}\label{sec:r to t}

In this section, we convert radial decay $1/\jr^p, p \le1$ in $\{ r<t/2\}$ to $1/\jv^p$. In particular, the fundamental solution to the wave equation gives a $1/\jr$ decay rate, which we can now convert to $1/\jv$. This holds for both $\p$ and vector fields of $\p$. 

\begin{lemma} \label{use of sled} We have
\begin{align*}
\|\p_{\le m}\|_{L^\iy(\inte)} 
	&\ls \f1{T^{3/2}} \|\jr \p_{\leq m+n}\|_{LE^1(\tinte)}.
\end{align*} 
\end{lemma}

\begin{proof} 
This estimate will follow as a consequence of \cref{Sob.emb} and from proving that 
\begin{equation}\label{eq:LE1bound}
\|\p_{\le m}\|_{LE^1(\tinte)}\ls \f1T \|\jr\p_{\le m+n}\|_{LE^1(\tinte)}.
\end{equation} 
The statement \cref{eq:LE1bound} hints at the fact that we will transfer (a limited amount of) $\jr$ decay into $\jv$ decay in the $LE^1$ local energy norm. From the $LE^1$ norm, we can recover pointwise bounds simply by explicit computation.%

\begin{remark}[It suffices to look at $\phi$ supported in $\inte$] \label{rem:support} In this proof we can assume that $\p$ is supported in $\inte$ because we can control the commutator $[P,\chi_{\inte}]$ adequately where $\chi_\inte$ is a cutoff function adapted to the region $\inte$. Henceforth, we will assume that $\phi$ is supported in $\inte$, although support in $\jr\le \xl T$ for any fixed $\xl>0$ would also be fine. 
\end{remark}

Let $m\geq 0$. Let $\xg_{(T,x)}(t')$ denote an integral curve of $S$, parametrized by unit speed, such that $t'=0$ corresponds to time $t = T$ and spatial position $x$. That is, it corresponds to the point $(T,x)$. 
By the fundamental theorem of calculus and Cauchy-Schwarz, we have
\begin{equation}
|\de \pmo(T,x)|^2
\ls \f1T \int_0^{T} 
|(\de \pmo)(\gamma_{(T,x)}(t'))|^2 + |(S\de \pmo)(\gamma_{(T,x)}(t'))|^2\,dt'.
\label{vfromsv}\end{equation}
(This bound clearly works for any smooth-enough function other than $\pmo$ as well.)
Next, integrating \cref{vfromsv} on $\{x : r< \xl t\}$ for some $\xl>0$, say, $\{ x : r \le t \}$, 
\begin{align*}
\int_{C_T\cap \{t=T\}} |\de \pmo(T,x)|^2\,dx 
	&\ls \f1{T} \iint_{C_T} |\de \pmo |^2 + |S \de \pmo |^2 \,dxdt \\
	&\ls \f1T \iint_{C_T} |\de\p_{\le m+2}|^2 \,dxdt.
\end{align*}
  
  A similar bound holds for $t = 2T$, where we now average over $[0,T]$ again but this time over the integral curves of $-S$, using $\xg_{(2T,x)}(t')$ as the argument for the function $\de\pmo$, with $t'=0$ corresponding to time $t = 2T$. Thus
$$|\de \pmo(2T,x)|^2
\ls \f1T \int_0^{T} |(\de \pmo)(\gamma_{(2T,x)}(t'))|^2 + |(S\de \pmo)(\gamma_{(2T,x)}(t'))|^2\,dt'.$$
Then, integrating over $\{x : r \le t\}$, we obtain the same upper bound $T\inv \|\de\pmt\|^2_{\lt(C_T)}.$

Hence by the solution $\phi$ satisfying \cref{wled for vf} (weak local energy decay for vector fields),
\begin{equation}\label{pre wwled est}
    \| \p_{\leq m}\|_{LE^1(\tinte)}  \ls \|\de\pmo(T)\|_\lt   \ls \f1{T^{1/2}} \| \de \p_{\leq m+2}\|_{L^2(C_T)}.
  \end{equation}

  \
  
Next, we bound $\|\de\p_{\le m+2}\|_{\lt(C_T)}$ using \cref{lem:wsled}.  Intuitively, \cref{lem:wsled} ``multiplies'' or ``boosts'' all integrands in \cref{wled for vf} by $\jr^{1/2}$. A first naive thought that comes to mind is to multiply the equation by $r\pa_r\phi$ to achieve this boost. This works, if we add a zeroth-order correction term $\phi$ to the multiplier. Unlike the unweighted multiplier, this weighted multiplier leads to unsigned constant-time boundary terms, hence we put both energy terms in the majorizer. %
\cref{lem:wsled} adds new information beyond \cref{wled for vf} only for sufficiently large values of $r$. 

We will sometimes use the notation $C_{T_1}^{T_2} := [T_1,T_2]\times \{x : r\le t\}$.

\begin{lemma}[$\jr^{1/2}$-weighted Weak Local Energy Decay for Vector Fields] \label{lem:wsled} Suppose that the solution to
$P\phi = f$ satisfies the weak local energy decay for vector fields, as proved in \cref{wled for vf}.  For all $0 \leq  T_1 \leq T_2$, we have
\begin{align}\label{wwled est}
\begin{split}
\| \de \phi_{\leq m}&\|_{\lt(C_{T_1}^{T_2})}
\ls \sum_{j=1}^2 \| \jr^{1/2} \de \pmo(T_j)\|_\lt  + \|\jr f_{\leq m+2}\|_{L^2[T_1,T_2]L^2}.
\end{split}
\end{align}

\end{lemma}

\begin{proof}We will take as assumptions those stated in part (1) of \cref{thm:main} and prove this result. This implies that this result also holds for part (2), because the assumptions in part (2) are stronger than those for part (1). 

\begin{itemize}
\item (The zero multiindex case)
We demonstrate the case $m=0$ first for simplicity. In this proof we shall need $\xs$ and $\xd$ to be strictly positive real numbers, as well as $\xd'>-1$, in contrast to the situation in \cref{dya bd}. 

We multiply $P\p=f$ by $r\pa_r\phi + \phi$ and integrate by parts in $[T_1,T_2]\times\R^3$. There is a number $q'>0$ such that
\begin{align}\label{computation}
\begin{split}
 \int |\de &\phi|^2 + O(\jr^{-q'})  |\de \phi|^2 + O(\jr^{-1-q'}) |\pao\p|^2 + O( \jr^{-2-q'}) |\phi|^2  \,dxdt \\
  &\ls \sum_{j=1}^2 \int_{\R^3} O(\la r \ra) |\de \phi(T_j,x)|^2 +  O(\la r \ra^{-1}) |\phi(T_j,x)|^2 \, dx + \int |rf\pa_r\phi| + |f\p| \,dxdt \\
  &\ls \sum_{j=1}^2 \int_{\R^3} O(\la r \ra) |\de \phi(T_j,x)|^2 \, dx + \int |rf\pa_r\phi| + |f\p| \,dxdt \\
\end{split}
\end{align}with the last statement following by a version of Hardy's inequality.

	For instance, with the term $\pa_\mu A^\mu$, 
$$\pa_\mu (A^\mu \phi) (r\pa_r \phi) \ls O(\jr^{-\xs}) |\de \phi|^2 + O(\jr^{-2-\xs})\phi^2.$$
This follows by combining the assumptions on $A,\pa A$ and $\pat A$ as stated in part (1) of \cref{thm:main}. Only an arbitrarily small $\xs>0$ is needed. 

Next, 
$$\iint |f| |r\pa_r \phi| = \iint |fr| \, |\pa_r \phi| \leq \|rf\|_\ltlt \|\pa_r \phi\|_\ltlt \le \f1\eps \|rf\|_\ltlt^2 + \eps\|\pa_r \phi\|_\ltlt^2$$ 
and we then bring $\eps\|\pa_r \phi\|_\ltlt$ onto the other side together with $\|\de \phi\|_\ltlt$. Similarly
$$\iint |f| |\phi| = \iint |rf| \, |\f{\phi}r| \le \f1\eps \|rf\|_\ltlt^2 + \eps\|\pa_r \phi\|_\ltlt^2.$$ 
We remark that it is possible to place $rf$ in $\lolt$ if we place $\pa_r\p$ in $L^\iy \lt$ (we can use Hardy's inequality for the zero order term). This alternate route leads to $\|rf\|_\lolt$ instead of $\|rf\|_\ltlt^2$ on the right-hand side. 

	For small $|x|$ values, our assumption of \cref{wled for vf} implies the bound on $$\|\de\p\|_{\lt[T_1,T_2]\lt(r \ls 1)}.$$
	
	By using the positivity of $q'$ on the left-hand side of \cref{computation} for large $|x|$ values, we can then obtain
\begin{equation}\label{sled m=0}
\begin{split}
\| \de \phi&\|_{\lt[T_1,T_2]\lt} \ls \sum_{j=1}^2 \| \jr^{1/2} \de \phi(T_j)\|_{L^2} + \min\lr{ \|\jr f\|_{L^1[T_1,T_2]L^2}^{1/2} ,   \|\jr f\|_{L^2[T_1,T_2]L^2} }.
\end{split}
\end{equation}
\cref{sled m=0} implies \eqref{wwled est} for $m=0$. 

\item (The higher multiindex case)
Next, we prove \eqref{sled m=0} but for $\phi_J, J\ne \vec 0$. 
By \cref{lem:comm0}, we have
\begin{align*}
P \phi_J 	&= f_{J} + O(\jr^{-1-q'}) \de \phi_{\leq |J|} + O(\jr^{-2-q'}) \phi_{\leq |J|-1}.
\end{align*}
We multiply this by $r\pa_r \phi_J + \phi_J$. Then we integrate in $[T_1,T_2] \times \R^3$. 
\begin{itemize}
\item
For small $r$, the estimate \eqref{wwled est} is implied by the weak local energy decay estimate for vector fields proved in \cref{wled for vf}, so to prove the desired conclusion \eqref{wwled est} it suffices to restrict attention to the case of large $r$. 

\item
For large $r$, owing to the positivity of $q'>0$, we may use the triangle inequality, the triangle inequality for integrals, Cauchy-Schwarz, and Hardy's inequality to absorb the terms
\[ 
\iint (r\pa_r \phi_J + \phi_J) \Big( O(\jr^{-1-q'}) \de \phi_{\leq |J|} + O(\jr^{-2-q'}) \phi_{< |J|} \Big) 
\]
into the left-hand side, namely into $$\| \de \phi_{\leq m}\|_{\lt([T_1,T_2]\times \{ r\le t\})}.$$
The positiveness of $q'$ provides the necessary smallness for the absorption. 

	We have explained how to take care of the extra terms arising from commutators in the higher multiindex case, namely the terms 
$O(\jr^{-1-q'}) \de \phi_{\leq |J|} + O(\jr^{-2-q'}) \phi_{\le |J|-1}.$
	For the remaining part of the equation, namely $P\phi_J = f_J + (\text{taken care of})$, we can just apply precisely the same procedure used to prove \cref{sled m=0} to $\phi_J$---that is, the first bullet point. Then we sum over $|J| \leq m$.
\end{itemize} 
\end{itemize}
\end{proof}

Applying \cref{lem:wsled} for $\p$, we have
$$\|\de \pmt\|_{\lt(C_T)} \ls \sum_{i =1}^2 \|\jr^{1/2} \de \p_{\le m+3}(i T)\|_\lt.$$
The next step will be to bound these weighted energy terms by $LE^1$ norms, picking up appropriate $T$ weights along the way.

By the fundamental theorem of calculus and Cauchy-Schwarz once more,
\[
\int \jr |\de \pmn(T)|^2\,dx
    \ls \f1{{T}^{1/2}} \int \jr^{1/2} |\de \pmn |^2 + \f1{T} \jr^{3/2} |S \de \pmn |^2 \,dxdt.
\]
By \cref{rem:support}, we assume that $\jr\ls T$, which lets us bound
\[
\f1{T^{3/2}}\int \jr^{3/2}|S\de\p_{\le m+n}|^2\,dxdt 
\ls \f1T \int \jr |S\de\pmn|^2 dxdt 
\ls \f1{T} \|\jr\p_{\le m+n+n'}\|_{LE^1(C_T)}^2 .
\]for some $n'$. %

We are not able to directly bound ${T}^{-1/4}\|\jr^{1/4} \de \pm\|_{L^2(C_T)}$ by ${T}^{-1/2} \|\jr \pmn\|_{LE^1(C_T)}$. %
Instead, we treat this term perturbatively for small $r$, and for a fixed finite number of large $R$ regions, where $r\sim R$, we can make this bound. Thus let us decompose
$$T^{-1/4}\|\jr^{1/4}\de\pm\|_{\lt(C_T)} = \sum_{R} T^{-1/4} \|\jr^{1/4}\de\pm\|_{\lt([T,2T]\times A_R)}.$$
When $R \ll T$ we absorb this term into the left hand side. For all values of $R$ with %
$R \sim T$, we are able to directly bound by ${T}^{-1/2} \|\jr \pmn\|_{LE^1(C_T)}$. Thus from
\begin{align*}
\|\de \p_{\le m+2}\|_{\lt(C_T)}
	&\ls {T}^{-1/4}\|\jr^{1/4} \de \p_{\le m+n}\|_{L^2(C_T)}
		+ {T}^{-1/2} \|\jr \p_{\le m+n}\|_{LE^1(C_T)} 
\end{align*}
we are able to conclude 
$$\|\de \p_{\le m+2}\|_{\lt(C_T)}\ls {T}^{-1/2} \|\jr \p_{\le m+n}\|_{LE^1(C_T)}$$
so that 
$$\|\p_{\le m}\|_{LE^1(\inte)} 
	\ls T^{-1} \|\jr \p_{\le m+4}\|_{LE^1(C_T)},$$ 
	which proves \cref{use of sled}. 
	Note that we could actually have obtained the upper bound in $LE^1(\tinte)$, rather than $LE^1(C_T)$, i.e. 
	$$\|\p_{\le m}\|_{LE^1(\inte)} 
	\ls T^{-1} \|\jr \p_{\le m+4}\|_{LE^1(\tinte)},$$
\end{proof}

\begin{theorem}\label{r-t} Let $\phi$ solve the main equation \eqref{main}.
If \begin{equation}\label{hypothesis r-t}
\phi_{\le M} \ls \jr^{-p}\jt^{-q}\nm^{-\eta}
\end{equation}
for some real $p\le 1$, $q,\eta \in \R$ and a (sufficiently large) fixed $M\in \N$.
 then 
 \begin{equation*}
 \phi\ls \jt^{-p-q}\nm^{-\eta}.
 \end{equation*}
\end{theorem}

\begin{proof}
For all $(t,r)$ pairs with $r$ sufficiently large relative to $t$, say $r > t/2$, the conclusion follows since $\jr \sim \jt$. 

For the other region, $\inte$, this follows from the proof of \cref{lem:r-t}, because in $\inte$, $\nm \sim \jt$. 

\begin{lemma}\label{lem:r-t}
Let $\phi$ solve the main equation \eqref{main}. 
If \begin{equation}
\phi_{\le M} \ls r^{-p}\jt^{-q}
\end{equation}
for some real $p,q \in \R$ and a (sufficiently large) fixed $M\in \N$ where
$p\le1$,
 then 
 \begin{equation*}
 \phi\ls \jt^{-p-q}.
 \end{equation*}
\end{lemma}

\noindent\textit{Proof of \cref{lem:r-t}.} %
We compute the norms involved on the right-hand side in \cref{use of sled}. 
	The rest of this proof works for not only $\inte$, which is the region we compute in, but actually in $[T,2T] \times \{ r \le \xl t\}$ for any fixed $\xl > 0$.
The right-hand side norm of \cref{use of sled} is
\begin{align*}
\|\jr \p_{\leq m+n}\|_{LE^1(\inte)}  
	&= \| \de(\jr \p_{\le m+n} )\|_{LE(\inte)} + \| \p_{\le m+n}\|_{LE(\inte)} \\
	&\ls \|\jr \de \p_{\le m+n}\|_{LE(\inte)} + \| \p_{\le m+n}\|_{LE(\inte)}  \\
	&\ls \| \p_{\le m+n}\|_{LE(\inte)}
\end{align*} 
where the last line is a consequence of \cref{derdecbtr} applied uniformly across the collection $\{\crt : 1 \le R < 3T/8 \}$ of dyadic regions.
Thus
$
\|\phi_{\le m}\|_{LE^1(\inte)}  
	\ls \f1T  \| \phi_{\leq m+n}\|_{LE(\inte)}  
$. 

Next, we bound $\|\phi_{\leq m+n}\|_{LE(\inte)}$ and finish the proof. We shall use pointwise bounds on $|\pm|$, and not just on $|\p|$, here: 
\begin{itemize}
\item
For $R>1$ we have \begin{align*}
\sup_{1 < R < 3T/8} &\left( \int_T^{2T} \int_R^{2R} \f1\jr \left( \jt^{-2q} r^{-2p} \right) r^2drdt \right)^{1/2} 
	\sim \sup \left(T^{-2q} \int_T^{2T} \int_R^{2R} \f1\jr \left( r^{-2p} \right) r^2drdt \right)^{1/2} \\
	&\ls  T^{1/2 - q} \sup_{1 < R < 3T/8} \f1{R^{p-1}} \ls T^{1/2-q} \f1{T^{p-1}} \text{ since }p\le1.
\end{align*}
\item
For $R=1$ we have 
\[
\lr{ \int_0^2 \f1{\jr^{2p-1}}dr }^{1/2} \ls_p 1
\]
for any $p\in \R$. 
\end{itemize}
Thus
\begin{align*}
\f1{T^{3/2}} \|\pmn\|_{LE(\inte)} 
 	&\ls \f1{T^{3/2}}  \f1{T^{ \min(0,p-1) + q - 1/2 }} \\
	&= \f1{T^{\min(1,p)+q}}
\end{align*}
hence $\|\p_{\le m}\|_{L^\iy(\inte)} \ls T^{-p-q}$ if $p \le1$. 
\end{proof}
This establishes the proof of \cref{r-t}.
\hfill $\Box$

\begin{corollary}
Let $k\ge1$ be an integer. If $\p$ solves $P\p=f$ and
$\p_{\le M} \ls \sum_{j=1}^k \jr^{-p_j}\jt^{-q_j}\nm^{-\eta_j}
$and the conditions on the exponents $p_j,q_j,\eta_j$ and $M$ in \cref{r-t} above are satisfied, 
then 
$\p \ls \sum_{j=1}^k \jt^{-p_j-q_j}\nm^{-\eta_j}.
$\end{corollary}

\begin{proof}The proof is a straightforward consequence of what has already been done. One can use elementary inequalities to handle sums instead of single summands in the computations above, and the estimates still hold. %
\end{proof}

\section{The upper bound in $\{ r<t\}$}\label{sec:inte}

We consider \eqref{wm eqn} with $r<t$. We now show the desired final decay rate in \cref{thm:main}, namely
$$\wm \ls \f1{\jr\nm^{1+\min(1+\xs,\xd,\xd')}}.$$ 

\begin{proposition}Assume that $r<t$. 
Assuming the hypotheses of part 2 of \cref{thm:main}, $$\wm \ls \f1{\jr\nm^{1+\min(1+\xs,\xd,\xd')}}.$$ 
Assuming the hypotheses of part 1 of \cref{thm:main}, $$\wm \ls \f1{\jr\nm^{1+\min(\xs,\xd,\xd')}}.$$ 
\end{proposition}

\begin{proof}
In $\{ \rho\geq s\}$, the argument has essentially been done in \cref{sec:exte}; when integrating in $\{\rho>s\}$, one plugs in the final pointwise decay rates for vector fields of $\p$ obtained in \cref{sec:exte}, only to get the final pointwise decay rates as output. %

In $\{ \rho<s\}$, we let $\nu := \min(1+\xs,\xd,\xd',1-) = \min(\xd,\xd',1-)$\footnote{we need $\nu<1$ because we will be using the fact that in $\rei$, we have \[\sum_{R\in\rei}\int_\co \rho|\bo \wm|dsd\rho \ls \f1{\nm^{\xb+\eta+\x-3}} \] (using the notation from \cref{conversion}).}. Note that $2+\nu <3$, allowing us to apply \cref{conversion}'s Region 1 bound $1/\nm^{\xb+\eta+\x-3}$ if we put $\jv$ as $\jt$ in \eqref{onea} when applying \cref{conversion}. For the rest of the proof, the strategy will be to improve by increments $\nu$ which are strictly less than 1. Below in the proof, we split into the cases where $\min(\xd,\xd')$ is either $<1$ or $\ge1$, but the main idea in either case is really the same, since in the latter case we simply introduce an artificial decrement $\ti\eps\ll1$ to make $\nu$, which equals $1-\ti\eps$ in that case, smaller than 1. 

By \cref{dya bd}, we have 
\begin{equation}\label{onea}
\Box \wm \ls \jr^{-2-\nu} \jv\inv \nm^{1/2}.
\end{equation}
By \cref{conversion}, 
 \[ \jr \wm \ls \nm^{1/2 - \nu }. \] 
We have gained $\nm^{-\nu}$. 
Hence by \cref{r-t},
 \[\Box \wm \ls \jr^{-2-\nu} \jv\inv \nm^{1/2 - \nu }, \] 
and this process can be continued as long as the uppermost case thresholds in the definition of $\xk$ are not met. Suppose $n'>0$ is an integer for which this threshold is not met; then after performing this procedure $n'$ times, \[\Box \wm \ls \jr^{-2-\nu} \jt\inv \nm^{1/2 - n'\nu }. \] %
Now we define $n'$ to be  
\[
n' := \max\{ n \in\N: 1/2 + n\nu  < 1\}. \]
There are two cases:

\begin{enumerate}
\item
If $$1 < 1/2 + (n'+1)\nu < 1+\nu$$ then write $$1/2  + (n'+1)\nu = 1 + \xl\nu$$ where $0<\xl<1$; thus 
	\[\wm\ls \jr\inv \nm^{-1-\xl\nu}.\] 
Then 
\begin{align*}
\Box \wm &\ls \jr^{-2-\nu} \jv\inv \nm^{-1-\xl\nu} \\
	&\le \min \{\jr^{-2-\nu} \jv\inv \nm^{-1-\xl\nu}, \jr^{-3-\nu} \nm^{-1-\xl\nu} \} \\
	&=: \min\{a,b\}. 
\end{align*}
 	\cref{conversion} implies that in Region $\rei$, we have the bound by 
\[ 1/\nm^{\xb + \eta + \x - 3}.\] We use $a$ to get the bound in $\rei$ by \[\nm^{-1-\nu-\xl\nu}. \] 
On the other hand, we use $b$, with $\x=3+\nu$ and $\xb+\eta = 1+\xl\nu$ to get a Region $\reii$ bound by \[
\nm^{-1-\nu}. \]
Thus
\begin{align*}
\jr \wm &\ls  1/\nm^{\xb+\eta+\x-3} +  \nm^{-(\xb+\x-2)} \xk(\eta,t-r) \\
	&= \nm^{-1-\nu-\xl\nu} + \nm^{-1-\nu}\\
	&\ls\nm^{-1-\nu}. 
\end{align*}
\item
If $$1/2  + (n'+1)\nu = 1,$$ thus $\wm \ls \jt\inv \nm\inv$, we have $$\Box \wm \ls \jr^{-2-\nu} \jt \inv \nm \inv.$$ Hence 
\begin{align*}
\jr \wm 
	&\ls \nm^{-(\xb+\eta)} \xk(\x-1,t-r) + \nm^{-(\x+\xb-2)} \xk(\eta,t-r) \\
	&= \nm^{-(\xb+\eta)} + \nm^{-(\x+\xb-2)} \log\nm \\
	&= \nm^{-2} + \nm^{-1-\nu} \log \nm \\
	&\le 2 \nm^{-1-\nu} \log\nm \\
	&\ls \nm^{-1-\xl\nu}
\end{align*} for any $0<\xl<1$, which now puts us in case (1).
\end{enumerate}  The proof is complete when $\min( \delta,\delta' ) < 1$. 

\

\noindent \textit{Part two of the proof: The case where $\min(\xd,\xd')\ge1$, that is, all three parameters $\xd,\xd'$ and $1+\xs$ are at least 1.}
Suppose that $\min(1+\xs,\xd,\xd')\ge1$.   We shall still work with an increment $\nu$ that is less than 1.  Rather than using ``$1-$'' in the definition of $\nu$, we write $\nu$ as the definite number $\nu := 1-\ti\eps$
where $\ti\eps>0$ is a small number.
Then 
\begin{align*}
\bo \wm 
	&\ls \jv\inv \nm^{-1-\nu}\min\{\jr^{-2-\nu},\jr^{-2-\min(1+\xs,\xd,\xd')} \} \\
	&= \jv\inv \nm^{-1-(1-\ti\eps)}\min\{\jr^{-2-(1-\ti\eps)},\jr^{-2-\min(1+\xs,\xd,\xd')} \} 
\end{align*}
where we wrote down the trivial minimum of the two powers of $\jr$ to emphasise the fact that we will be using $\x=2+\min(1+\xs,\xd,\xd')$ for $\reii$ but $\x=2+\nu$ for $\rei$. 
Thus by \cref{conversion},
\begin{align*}
\jr \wm 
	&\ls \nm^{-(\xb+\eta+\x-3)} + \f{\xk(\eta,t-r)}{\nm^{\x+\xb-2}}	\\
	&\ls \nm^{-(\xb+\eta+\x-3)} + \f1{\nm^{\x+\xb-2}}	\\
	&= \nm^{-(\xb+\eta+\x-3)} + \f1{\nm^{ (2+\min(1+\xs,\xd,\xd') )+ (1) - 2 }} \\
	&= \nm^{-1-2\nu}	 + 	\f1{\nm^{1 + \min(1+\xs,\xd,\xd')}}. 
\end{align*}
It remains to prove the desired bound in Region $\rei$, and it is safe to ignore the $\reii$ portion of the bound henceforth because the $\xb$ and $\x$ exponent components of $\bo \wm$ remain stable while $\eta>1$ will stay larger than 1, and in $\reii$ we use the bound $\xk(\eta,t-r)/\nm^{\x+\xb-2}$. We note that no more improvement is possible in $\reii$ using  \cref{conversion}. 

In $\rei$, this iteration continues until \[\wm\ls\jr\inv\nm^{-1-n''\nu} \]
where 
\[
n'':=\max\{n\in\N:n(1-\ti\eps) < \min(1+\xs,\xd,\xd') \}, \] e.g., $n''=1$ if the two numbers $\min(1+\xs,\xd,\xd')$ and $1-\ti\eps$ are both close to 1. 

One way to view this situation is that there are two cases:
\begin{enumerate}
\item
If 
\[(n''+1)(1-\ti\eps) > \min(1+\xs,\xd,\xd') \]
then we obtain the bound $\nm^{-1-(n''+1)\nu} = \nm^{-1-(n''+1)(1-\ti\eps)}$ in $\rei$ by using
\begin{align*}
\f1{\nm^{\xb+\eta+\x-3}} &= \f1{\nm^{1 +(n''+1)\nu}}\\
	&\le \f1{\nm^{1+\min(1+\xs,\xd,\xd')}}.
\end{align*}
\item
If 
\[(n''+1)(1-\ti\eps) = \min(1+\xs,\xd,\xd') \]
then we obtain the final display in item (1) but with equality rather than inequality, and we halt. 
\end{enumerate}
This completes the proof for $\wm$ when $r<t$. 
\end{proof}

\begin{remark}[Lockstep] \label{lockstep}
	If $\nu :=\min(\xs,\xd,\xd',1-)$ then essentially an identical proof follows for proving $\wm\ls \f1{\jr\nm^{1+\min(\xs,\xd,\xd')}}.$ The case partition is then (a) part one: $\min(\xs,\xd,\xd')<1$, (b) part two: $\min(\xs,\xd,\xd')\ge1$. Everything else follows when one replaces $1+\xs$ in the appropriate locations in the proof above by $\xs$. 
\end{remark}

\section{Cone bounds} \label{sec:cone bounds}
In this section we show how we prove the final decay rate in the main theorem for the terms involving the metric coefficients $\hab$ that are supported near the cone $\{r=t\}$.

Recall that we write $\ti B\pa\pm = \pa(\ti B \pm) - (\pa\ti B)\pm$. Let $j=0$ (respectively $j=1$) correspond to the hypotheses of part 1 (respectively part 2) of \cref{thm:main}. Near the cone, we rewrite \eqref{pr} as
\begin{align}\label{rewrite}
\Box\pm &= ( |\ti V| + |\pa \ti B| + |\ti{g^\xo}|)|\pmn| + \pat( \cone (\ti h \pat+\ti A+\ti B)\pm), \\
(\pm(0),\bar N\pm(0)) &= (0,0).
\end{align} 
and use \eqref{D decomp}. 
Note that $|\ti V| + |\pa \ti B| + |\ti{g^\xo}| \ls |\ti V| + |\pat\ti B| + \jr\inv|\ti B_{\le n}| + |\ti{g^\xo}|= O( 1/\jr^{2+\min(\xs+(j-1),\xd,\xd')} )$ assuming the hypotheses of part $j$ of \cref{thm:main}, $j=1,2$. 

It suffices to prove pointwise decay estimates for
$$\bo v_{(m,1)} = \cone \ti h \pat\pm, \qquad \bo v_{(m,2)} = \cone (\ti A +\ti B)\pm,$$

Let $\ti v\in\{ v_{(m,j)}  : j=1,2\}$.  We now prove 
\begin{proposition}We have
$$\pat\ti v \ls \f1{\jr\nm^{1+\min(\xs,\xd,\xd')}}$$
under the assumptions of part 1 of the main theorem.
\end{proposition}

\begin{proof}
If $\chi := \cone$ and $f\in\{\chi (\ti A+\ti B) \pm,\chi\ti h\pat\pm\}$, then by \cref{derdecbtr} and assumptions on $h$ and $A$ we have 
$$|f(s,\rho)| +|Sf| + |\no \pa_\rho f| \ls \f1{\jrho^{1+\xs}} |\pmn|$$ %
	The iteration for $\ti v$ is as follows. Note that $\supp\chi \subset \{|s-\rho| \ls \nm\}$. One simplifying observation %
	is that $\rho \ge c \nm$ in all $f$ cases with $c \ge 1/4$; in $r<t$, $\supp\cone$ for instance, $\rho \ge |t-r|/4$, which has smallest $c$ value amongst all cases (for example, if $r>t$ then $c=1$ and if $r<t$ then in $\supp\cone$, $c = 1/2$). This and the fact that the horizontal (i.e. $\rho$) diameter of $\supp\chi$ is $O(\nm)$ leads to simpler integrations in $\rho$. 
	
We begin with the bound \eqref{u/v decay} for the functions $\pmn$. By \cref{derdecbtr} (used to handle terms that have the operator $\on \pa_\rho$ in the integration inside $\dtr$) and \cref{conversion}, 
$$ \pat \ti v \ls \f{\mn^{1/2-\xs}}{\jr}.$$
Thus we run the iteration %
 with exponent $a:=\min(\xs,\xd,\xd',1-)$---see \cref{lockstep}. 
By \cref{conversion}, after $N$ steps one gets $$\jr\pat\ti v\ls \f1{\nm^{1+\xs}}\nm^{3/2-\xt}$$ where $\xt=Na<3/2$ is the gain at the $N$-th step of the lockstep. The procedure is similar to $\wm$'s case, and in the end we get $$\pat\ti v \ls \f1{\jr\nm^{1+\min(\xs,\xd.\xd')}}.$$

We have $\pa(\ti A \pm) + \ti B \pa\pm = \pa( [\ti A + \ti B] \pm) - (\pa \ti B) \pm$.
For $\bar v$ solving 
$$\Box \bar v = \chi \ti A\pm,$$
the bound on $\pat \bar v$ is just 
an argument that is an application of \cref{conversion} similar to what has been done. 
	We write $B\pa\p = \pa(B\p) - (\pa B)\p$; then the arguments already shown (along with the assumptions on $B$ in \cref{thm:main}) give the bound on $\pat W$ for $\Box W=\chi \ti B \pm$. 
	This concludes the proof.
\end{proof}

\begin{proposition}
We have
$$\pat\ti v \ls \f1{\jr\nm^{1+\min(1+\xs,\xd,\xd')}}$$
under the hypotheses of part 2 of the main theorem.
\end{proposition}
\begin{proof}
We now prove $$\pat\ti v \ls \f1{\jr\nm^{1+\min(1+\xs,\xd,\xd')}}$$ assuming more on the time derivatives of our coefficients and also a little more on $A$ and $B$; see part (2) of \cref{thm:main}. 
	For the first-order terms, we again write $\pa(\ti A\pm) + \ti B \pa\pm = \pa( [\ti A + \ti B] \pm) - (\pa \ti B) \pm$; since $A$ and $B$, respectively $\pat A$ and $\pat B$, belong to the same $S^Z$ class, with one higher rate of $\jr$ decay relative to the hypotheses of part 1 of the main theorem, we are done by the previous proof and we henceforth focus on the metric coefficients.
By the product rule,
\begin{align}\label{rewriting}
\pat(\chi \ti h \pat \pm) 
	&= \pat^2(\chi \ti h \pm) - \pat( \pat(\chi \ti h) \pm)\end{align}
where $\pat(\chi \ti h) = O(\jr^{-2-\xs})$ since $\pat h\in S^Z_\text{cone}(\jr^{-2-\xs})$. 
For $\Box U = -\pat(\pat(\chi\ti h)\pm)$, the pointwise decay rates for $U$ follow from techniques already shown.

For $\Box \pat^2 u=\pat^2(\chi \ti h \pm)$, we have
\begin{align*}
\mn u_{tt} &\ls |Lu_t| + |Su_t| \ls |\pat (Lu) |+ |\pat u| +| \pat(Su)| + \f1\jr \sum_{\bar Z\in\{\xO,S\}} |\bar Zu|
\end{align*}where $L$ denotes the Lorentz boosts.
It suffices to bound the first three terms on the right hand side by $$O(\f{\nm^{3/2-\xt-q}}\jr)$$ if $$h\in S^Z(\jr^{-q}).$$

We now bound $|\pat u|$. 
Writing $\Box u=\chi \ti h\pm$ with simplified notation henceforth as $\chi h\p$ or $\chi h\pm$,
\begin{align}\label{c1}
\begin{split}
\jr\nm u_t &\ls \jr(|Su|+|Lu|) \\
	&\ls \int_\dtr |\chi h \p| + |S(\chi h\p)| + \on |\pa_\rho(\chi h\p)| \ \rho dsd\rho\\
	&\ls \f1{\nm^q}\nm^{5/2 -\xt}, \ \ q = 1+\sigma
\end{split}
\end{align}
where the last line follows from 
 \cref{conversion} and \cref{derdecbtr}. 

The same calculation shows that $(Su)_t$ is also bounded by this, since by replacing $u$ by $Su$ above we still find the same upper bounds for the integrand; this is because the three functions $S^j(\chi h\p), j=0,1,2$ satisfy the same bounds, and the analogous integral is 
\begin{align}\label{c2}
\begin{split}
\jr\nm (Su)_t &\ls \int_\dtr |\Box LSu| + \lr{|\Box S^2u|} \ \rho dsd\rho \\
	&\ls \int \on \sum_{k=0}^1 |\pa_\rho S^k(\chi h\p)| +	\lr{\sum_{j=0}^2 |S^j(\chi h\p)| } \ \rho dsd\rho \\
	&\ls \f1{\nm^q}\nm^{5/2 -\xt}.
\end{split}
\end{align}

The function $(Lu)_t$ also obeys the same bounds as $u_t$, and the analogous integral is 
\begin{equation}
\label{c3pre}
\jr\nm(Lu)_t  \ls \int_\dtr |\Box SLu| + |\Box L Lu| \ \rho dsd\rho.
\end{equation}
We have
\begin{align}\label{c3}
\begin{split}
\int_\dtr |\Box L Lu| \ \rho dsd\rho
	&\le \int  \lr{2|\Box Lu|} + | LL\Box u| \ \rho dsd\rho \\
	&\ls \int \lr{|S(\chi h\p)| + \on|\pa_\rho(\chi h\p)|} + |LL\Box u| \ \rho dsd\rho\\
	&\ls \int \lr{|S(\chi h\p)| + \on|\pa_\rho(\chi h\p)|} \\
		&+ |S^2(\chi h\p)| + |S(\on\pa_\rho(\chi h\p))| + \on|\pa_\rho S(\chi h\p)| \\
		&+ \on|\pa_\rho(\on\pa_\rho (\chi h\p))| \ \rho dsd\rho\\
	&\ls \f1{\nm^q}\nm^{5/2-\xt}, \ \ q = 1+\sigma
\end{split}
\end{align}where we changed $\pa_\rho$ to $\pa_s$ by \eqref{D decomp} and used the assumption on $\pat^2h$. 
The final line follows by \cref{conversion}. The final integrand term (and specifically, when the two derivatives both fall on $h$) is 
the sole instance the extra assumption on $\pat^2 h$ in \cref{thm:main} is used. %

For $\Box SLu$ we have
\begin{align*}
\Box SLu &= S\Box Lu  + 2\Box Lu \\
	&= SL\Box u + 2\Box Lu \\
	&= LS\Box u + O(t\pa(\chi h\p)) +2\Box Lu 
\end{align*}
where $O(t\pa(\chi h\phi))$ arises from $[S,L]$ and can be broken into three cases: this function takes one of the three forms $t\pai(\chi h\p) = \pai(t\chi h\p), x_i\pat(\chi h\p) = \pat(x_i\chi h\p)$, and $t\f{x_i}r \pa_r(\chi h\p) = \pa_r(t \f{x_i}r \chi h\phi)$. We may then replace $\pai,\pa_r$ in the first and third cases by $\pat$ via \eqref{D decomp}. Then all these terms on the right hand side yield the upper bound $\nm^{5/2-\xt-(1+\xs)}$ via \cref{conversion}; to see this, it suffices to consider a solution of $\Box w = t \chi h \phi$ and prove bounds for $\pat w$. 

For $LS\Box u$ on the other hand, 
\begin{align*}
\int_\dtr |S^2(\chi h\p)| + \on |\pa_\rho S(\chi h\p)| \ \rho dsd\rho \ls \f1{\nm^q}\nm^{5/2-\xt}, \ \ q=1+\xs
\end{align*}also. For $\Box Lu$, this upper bound was proved earlier. Thus
$$\jr\nm(Lu)_t \ls \f1{\nm^q}\nm^{5/2-\xt}.$$

In summary, $(Lu)_t$ and $(Su)_t$ obey the same bound as $u_t$, because $\Box (Lu)_t$ and $\Box (Su)_t$ obey the same bounds as $\Box u_t$ and the claim then follows from \cref{conversion}. Thus
$u_{tt} \ls \f1\jr \nm^{1/2-\xt-q}, \ \ q = 1+\xs$ and the iteration finishes with $\pat \ti v \ls \f1{\jr\nm^{1+\min(1+\xs,\xd,\xd')}}$. 
Notice that this argument works for any positive value of $q$. 
\end{proof}

\section*{Acknowledgements}
I am grateful to Mihai Tohaneanu for suggesting this problem and for discussions related to this work, and to Katrina Morgan and Jared Wunsch for discussing their paper \cite{MW}.

\end{document}